# SURVIVAL AND COMPLETE CONVERGENCE FOR A SPATIAL BRANCHING SYSTEM WITH LOCAL REGULATION

By Matthias Birkner and Andrej Depperschmidt

*Weierstrass Institute for Applied Analysis and Stochastics and Technische Universität Berlin*

We study a discrete time spatial branching system on $\mathbb{Z}^d$ with logistic-type local regulation at each deme depending on a weighted average of the population in neighboring demes. We show that the system survives for all time with positive probability if the competition term is small enough. For a restricted set of parameter values, we also obtain uniqueness of the nontrivial equilibrium and complete convergence, as well as long-term coexistence in a related two-type model. Along the way we classify the equilibria and their domain of attraction for the corresponding deterministic coupled map lattice on $\mathbb{Z}^d$.

**1. Introduction and main results.** An important problem from the field of mathematical ecological modeling is to find plausible stochastic models on the level of individuals for the time evolution of a "population," say, of animals or plants, which live, move—in the case of plants, we think rather of the dispersal of seeds—and reproduce in a 2-dimensional space, subject to individual random fluctuations. The mathematically simplest class of stochastic models one might come up with, namely, branching random walk and its relatives in which individuals do not interact, are not adequate because in dimension 2, they virtually never exhibit stable long-time behavior: it is well known (see, e.g., [20] and the discussion there) that they will die out locally if the branching is (sub-)critical, and grow locally beyond all bounds if it is supercritical.

To describe an "old" population, which corresponds mathematically to a nontrivial equilibrium situation, one has to introduce some interactions among individuals, which is of course also natural from the modeling perspective. A very drastic solution that is frequently used in so-called stepping stone models (see, e.g., [19]), in the context of population genetics models, is to force the population size, or the population size per deme in a spatially extended scenario, to be constant, that is, each birth is exactly matched by a death in the population. A number of ecological models have been introduced and studied rigorously in the context of interacting particle systems and probabilistic cellular automata on $\mathbb{Z}^d$ (see, e.g., [7, 8, 18]). In these models the state of the system at discrete or continuous time $t$ is described by a function $\xi_t : \mathbb{Z}^d \to S$, where in most cases $S$ is some finite set of possible types. The interpretation is that a site $x$ is vacant at time $t$ if $\xi_t(x) = 0$ and occupied by *one* individual of type $i$ if $\xi_t(x) = i$ for some $i \in S \setminus \{0\}$. An individual at site $x$ changes its type or dies at a certain rate (with certain probability in discrete time setting), which depends on the neighborhood of $x$. Reinterpreting the type as occupancy numbers, this class can, in principle, accommodate models with a fixed a priori upper bound on particle density. On the other hand, it seems more natural to allow arbitrary population sizes or densities, and introduce a self-regulation mechanism which, for example, makes individual reproduction super-critical in presently sparsely populated regions and subcritical in crowded areas—accounting for stress or competition for resources. Such models with explicit space have been studied in the ecological literature (see, e.g., [2, 13]), mostly using computer simulations and heuristic arguments (see, e.g., [3, 17] for a comparison of different approaches). Recently, some variants of models of locally regulated populations have been studied in the mathematics literature and the possibility of long-time survival in certain parts of the parameter space has been rigorously proved for a continuous mass model [1, 10, 11].

We add to this literature a variant where particles live in discrete demes (arranged on $\mathbb{Z}^d$) in nonoverlapping generations, which looks as follows: In the absence of competition, an individual has on average $m > 1$ offspring. Due to competition, for example, for local resources, the average reproductive success of an individual at position $x$ is reduced by an amount of $\lambda_{xy} \geq 0$ by each individual at position $y$. Here $\lambda_{xy}$ is a finite range kernel on $\mathbb{Z}^d$. Thus, an individual at $x$ in generation $n$ will have a random number of offspring with mean given by

$$(1) \qquad \left(m - \sum_{y \in \mathbb{Z}^d} \lambda_{xy} \xi_n(y)\right)^+,$$

where $\xi_n(y)$ denotes the number of individuals at spatial position $y$ in generation $n$, and for $r \in \mathbb{R}$, we write $r^+ = \max\{0, r\}$. In particular, if the occupancy of neighboring sites is so high that the term in brackets is negative, no offspring are generated at site $x$ in this generation. For definiteness and simplicity, we assume that the actual number of offspring, given the present configuration, is Poisson-distributed with the above mean, and independent for different individuals. Once created, offspring take an independent random walk step according to a kernel $p$ from the location of their mother. In this way, our model incorporates individual-based random fluctuations in the number and spatial dispersal of offspring.

A formal specification of the model is given as follows: We assume that the motion/dispersal kernel $p = (p_{xy})_{x,y \in \mathbb{Z}^d}$ and the competition kernel $\lambda = (\lambda_{xy})_{x,y \in \mathbb{Z}^d}$ satisfy the following conditions:

(A1) The kernel $(p_{xy})_{x,y \in \mathbb{Z}^d} = (p_{y-x})_{x,y \in \mathbb{Z}^d}$ is a zero mean aperiodic stochastic kernel with finite range $R_p \geq 1$, that is, for all $x, y \in \mathbb{Z}^d$: $p_{xy} = 0$ for $\|x - y\|_\infty > R_p$.
(A2) $0 \leq \lambda_{xy} = \lambda_{0,y-x}$, $\lambda_0 := \lambda_{00} > 0$ and $\lambda_{xy} = 0$ for $\|y - x\|_\infty > R_\lambda$, where $1 \leq R_\lambda < \infty$.

For a configuration $\eta \in \mathbb{R}_+^{\mathbb{Z}^d}$ and $x \in \mathbb{Z}^d$, define

$$f(x; \eta) := \eta(x) \left( m - \lambda_0 \eta(x) - \sum_{z \neq x} \lambda_{xz} \eta(z) \right)^+ \tag{2}$$

and

$$F(x; \eta) := \sum_{y \in \mathbb{Z}^d} f(y; \eta) p_{yx}, \tag{3}$$

that is, $f(y; \eta)$ is the expected number of offspring generated at site $y$ and, thus, $F(x; \eta)$ is the expected number of individuals at $x$ in the daughter generation if the present configuration is $\eta$. Let $N^{(x,n)}$, $(x, n) \in \mathbb{Z}^d \times \mathbb{Z}_+$ be independent standard Poisson processes on $\mathbb{R}_+$. Given $\xi_n$, the configuration of the $n$th generation, $\xi_{n+1}$, arises as

$$\xi_{n+1}(x) = N^{(x,n)}(F(x; \xi_n)), \qquad x \in \mathbb{Z}^d. \tag{4}$$

By well-known properties of the Poisson distribution, this definition is consistent with the intuitive description given above. Note that, technically, this model is a "probabilistic cellular automaton" (see, e.g., [5]) with countably infinitely many possible states at each site.

As for all $\eta \in \mathbb{R}_+^{\mathbb{Z}^d}$, we have $f(x; \eta) \leq m\eta(x)$, for $m \leq 1$, one can easily construct a coupling of $(\xi_n)$ with a subcritical branching random walk. In that case $(\xi_n)$ becomes extinct in finite time with probability 1 starting from any finite initial condition. Our first result states in the case $m \in (1, 4)$ that if the competition is weak enough, the population, starting from any nontrivial initial condition, will survive for all time with positive probability.

THEOREM 1. *For each $m \in (1, 4)$ and $p$ satisfying* (A1), *there are choices of positive numbers $\lambda_0^* = \lambda_0^*(m, p)$ and $\kappa^* = \kappa^*(m, p)$ such that if $\lambda_0 \leq \lambda_0^*$ and* $\sum_{x \neq 0} \lambda_{0x} \leq \kappa^* \lambda_0$, *then the population survives with positive probability, that is,*

$$\mathbb{P}_{\xi_0}[\forall n \in \mathbb{N}, \exists x \in \mathbb{Z}^d : \xi_n(x) > 0] > 0$$

*for all $\xi_0$ with $f(x; \xi_0) > 0$ for some $x \in \mathbb{Z}^d$. Furthermore, conditioned on nonextinction*

$$\liminf_{N \to \infty} \frac{1}{N} \sum_{n=1}^{N} \mathbb{1}_{\{\xi_n(0) > 0\}} > 0 \qquad a.s.,$$

*in particular, the origin (and, in fact, any site $x \in \mathbb{Z}^d$) will be occupied at arbitrarily large times.*

Note that this result as well as Theorem 3 and Corollary 4 below work in any dimension $d \geq 1$ (with threshold values $\lambda_0^*, \kappa^*$ depending on $d$), in particular, it establishes the possibility of long-term survival in $d = 2$.

The small competition coefficients mean that the system will typically be able to maintain a high number of particles per site. In this sense, our result concerns a "high density regime." Technically, we follow the natural path of a block construction in conjunction with comparison with oriented percolation, that might be paraphrased as "life plus good randomness leads to more life, so show that bad randomness has small probability." We call a space-time point *occupied* if there are enough particles there and not too many in the neighborhood (see Definition 6 for details). The definition is such that in the corresponding deterministic model (which is a "coupled map lattice" in dynamical systems jargon; see, e.g., [4] for a recent survey of this field)

(5) $$\zeta_{n+1}(x) = F(x; \zeta_n), \qquad x \in \mathbb{Z}^d, n = 0, 1, \ldots,$$

in which the Poisson variables are replaced by their means, an occupied site would after finitely many steps "colonize" its neighbors, that is, make them occupied as well. Then we control the probability that this remains the case under stochastic perturbation. Choosing small competition coefficients, we increase the "typical number of particles" per site in the deterministic model. Then we use the fact that the relative deviation of a Poisson random variable from its mean is typically small if the parameter is large. Finally, the finite range of competition and motion kernels allows to compare the set of occupied space-time sites with finite-range dependent oriented percolation on a suitable sub-grid of the space-time lattice.

The method can be adapted to a situation of two competing species to show that if in addition to the conditions of Theorem 1 the interspecific competition is weak enough, then long term coexistence is possible (see Theorem 8).

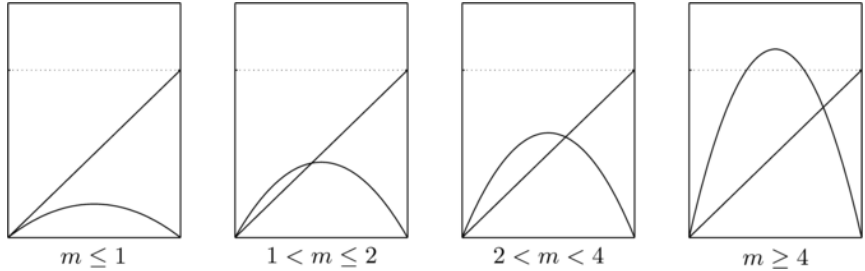

FIG. 1. *Function $\phi$ for different values of $m$.*

The logistic map $\phi(x) = x(m - \lambda x)^+$ and especially the one dimensional deterministic dynamical system

(6) $$x_{n+1} = \phi(x_n)$$

play an important role throughout the paper (see Figure 1 for a sketch of $\phi$ for various values of $m$). For example, in Theorem 1 the restriction to $m < 4$ comes from the fact that otherwise the function $\phi$ would not map the set $\{x \in \mathbb{R} : \phi(x) > 0\}$ into itself. The function $\phi$ has two fixed points, namely, 0 and $(m-1)/\lambda$. For $m \in (1,3)$, it is well known that 0 is repelling and $(m-1)/\lambda$ is attracting, that is, if $x_1 \neq 0$, then the sequence $(x_n)$ converges to $(m-1)/\lambda$, whereas for $m \geq 3$, there are no stable fixed points. It is interesting on its own that the former fact can be generalized to the coupled map lattice (5) which is a spatially extended version of (6): If the competition coefficients are small enough, the fixed point $\eta \equiv (m-1)/\sum_x \lambda_{0x}$ of the function $F$ is globally attracting in the sense of locally uniform convergence.

THEOREM 2. *Let $m \in (1,3)$, $p$, $\lambda$ satisfying* (A1) *and* (A2) *be given. Then there exists a positive number $\kappa^* = \kappa^*(m,p)$ such that if $\sum_{x \neq 0} \lambda_{0x} \leq \kappa^* \lambda_0$ and $f(x;\zeta_0) > 0$ for some $x \in \mathbb{Z}^d$, then $(\zeta_n)$ converges locally (i.e., pointwise w.r.t. $z \in \mathbb{Z}^d$) to $(m-1)/\sum_x \lambda_{0x}$.*

Notice that, under the assumptions of Theorem 2, we obtain a complete classification of the equilibria of (5) and their domains of attraction: If $(\zeta_n)$ does not hit the all zero configuration $\mathbf{0} \in \mathbb{Z}_+^{\mathbb{Z}^d}$ after the first step, it is attracted by $\eta \equiv (m-1)/\sum_x \lambda_{0x}$.

Obviously $\mathbf{0} \in \mathbb{Z}_+^{\mathbb{Z}^d}$ is an absorbing state for $(\xi_n)$, so the Dirac measure in this state is an invariant distribution for $(\xi_n)$. In view of Theorem 1, it is natural to ask if there exist nontrivial stationary distributions, and one might expect that if the process does not go extinct, its distribution converges to some unique invariant distribution. A powerful method to address this problem is coupling. Let $(\xi_n^{(1)})$ and $(\xi_n^{(2)})$ be versions of the process $(\xi_n)$ introduced in (4). Let $N_0^{(x,n)}$, $N_1^{(x,n)}$ and $N_2^{(x,n)}$, $(x,n) \in \mathbb{Z}^d \times \mathbb{Z}_+$, be independent standard Poisson processes. We define the coupling of $(\xi_n^{(1)})$ and $(\xi_n^{(2)})$ as follows:

$$
\begin{aligned}
\xi_{n+1}^{(1)}(x) &= N_0^{(x,n+1)}(F(x;\xi_n^{(1)}) \wedge F(x;\xi_n^{(2)})) \\
&\quad + N_1^{(x,n+1)}((F(x;\xi_n^{(1)}) - F(x;\xi_n^{(2)}))^+), \\
\xi_{n+1}^{(2)}(x) &= N_0^{(x,n+1)}(F(x;\xi_n^{(1)}) \wedge F(x;\xi_n^{(2)})) \\
&\quad + N_2^{(x,n+1)}((F(x;\xi_n^{(2)}) - F(x;\xi_n^{(1)}))^+).
\end{aligned}
\tag{7}
$$

THEOREM 3. *Let $m \in (1,3)$, and $p$, $\lambda$ as in (A1), (A2) be given. There are $\lambda_0^{**} = \lambda_0^{**}(m,p) > 0$ and $\kappa^{**} = \kappa^{**}(m,p) > 0$ such that if $\lambda_0 \leq \lambda_0^{**}$ and $\sum_{x \neq 0} \lambda_{0x} \leq \kappa^{**} \lambda_0$, then, conditioned on nonextinction of both populations, the coupling of $(\xi_n^{(1)})$ and $(\xi_n^{(2)})$ is successful in the sense that for each finite $\Lambda \subset \mathbb{Z}^d$, there is a random time $T$, such that*

$$\xi_n^{(1)}(x) = \xi_n^{(2)}(x) \qquad \text{for all } x \in \Lambda \text{ and } n \geq T.$$

Obviously we have $\lambda_0^{**} \leq \lambda_0^*$, $\kappa^{**} \leq \kappa^*$. We do not know if in the case $m \in (1,3)$ the inequalities are strict (but certainly the bounds obtained in the proof of Theorem 3 are much smaller than those obtained in the proof of Theorem 1). The proof of Theorem 3 again uses the block construction argument together with Lemma 14 which essentially follows from [6] where it was used in a similar spirit to prove uniqueness of the nontrivial invariant distribution for multicolor systems.

COROLLARY 4. *Under the conditions of Theorem 3 the process $(\xi_n)$ has two extremal invariant distributions. These distributions are translation invariant. Conditioned on nonextinction, $(\xi_n)$ converges in distribution in the vague topology to a random measure distributed according to the nontrivial extremal invariant distribution, that is, we have complete convergence.*

REMARK 5. (i) To our knowledge, we present here the first rigorous result showing the possibility of long-time survival in a locally regulated population in $d = 2$ for a particle-based model allowing multiple occupancy (but for particular cases in a continuous-time version, cf. [11], Proposition 6.4, where the competition acts strictly within-deme, and Proposition 7.9, where competition and dispersal kernel must be identical). Notice that because of the nonlocality of the competition and the discreteness of time the system considered here is not monotone (see [14], Chapter II.2, for background). As the mathematically rigorous investigation of spatial stochastic systems with local regulation terms is still in its infancy, we think it is justified to study the phenomenon in several mathematical guises. Furthermore, many species do live in discrete generations, and it is well known that discrete time dynamics can have a much richer behavior than their continuous time analogues. This shows up in our model as well; see point 4 below.

Being honest, one has to admit that the results of this paper, as well as those in [1, 11] are still too weak to capture many ecologically interesting phenomena. Up to now, all the rigorous results are more of a conceptual nature, showing that survival respectively coexistence of several types is possible if the interaction terms are weak enough, but giving little clues about what realistic sizes of threshold values enabling/excluding survival or coexistence might be. This stems from the fact that in order to apply comparison with finite-range dependent directed percolation, one usually has to keep far away from the true critical values. For example, we have little rigorous information about properties of the nontrivial equilibrium guaranteed by Corollary 4, apart from the fact that its mean is close to the deterministic prediction $(m-1)/\sum_x \lambda_{0x}$ when the competition terms are small. One would suspect that correlations decay exponentially, but we have no rigorous proof.

Thus, the contribution of these mathematical investigations to the question, how a population or several populations arrange themselves in space in order to survive in a (ecologically very interesting) situation of scarce resources and, hence, appreciable competition is at present rather limited. It appears that more powerful mathematical tools need to be invented in order to make rigorous progress in this direction.

(ii) The Poisson offspring distribution in our model is a somewhat artificial choice, which helps to streamline calculations, but is not essential for the result. To formulate a more general form of the model, one would need a one-parameter family of probability distributions (say, indexed by their mean) which includes sub- and supercritical distributions. A natural way would be to start with a fixed supercritical offspring distribution and then superimpose a "thinning" according to the local weighted density. A nice feature of the Poisson distribution is that we can in fact think of it in this way. Another feature of the Poisson distribution is that the variance of the total number of offspring produced at some site $x$ (given the present configuration) and its mean are the same. While it is natural for a "branching model" to assume that conditional variance and mean of the size of the new generation are of the same order, a general class of offspring distributions would allow for different proportionality factors.

(iii) Our results require that $\lambda_0$, the on-site competition coefficient, is (substantially) larger than the total competition with neighboring sites. Thus, they apply to a situation where most of the competition is felt by individuals within the same "colony." One can think, for example, of colonies arranged on $\mathbb{Z}^d$ and $\lambda_0$ governing a rather strong population regulation inside each colony, whereas the competition $\lambda_{0x}$, $x \neq 0$, with surrounding colonies is of a lower order.

This is certainly a technical condition which is not necessary for survival, but which intuitively helps quite a bit because it prevents the occupancy of a site from becoming so big that it would "eradicate" its neighborhood in the next step. Note that no such condition is necessary for the continuous-time continuous-mass result in Theorem 1.5, 2 b) in [10] (on the other hand, unlike [10], we do not need the requirement that the range of $\lambda$ must not exceed that of $p$).

Simulations suggest that the system may survive also when $\lambda_0$ and $\lambda_{0x}$, $0 < \|x\| \leq R_\lambda$, are the same or similar (but sufficiently small), but occupancy numbers will fluctuate much more wildly than in the scenario treated in Theorem 1. On the other hand, with a highly asymmetric competition kernel, one observes in simulations the appearance of "fronts" of occupied sites moving in the direction of smaller $\lambda$. This might indicate local extinction despite global survival when starting from a finite initial population in such a case.

(iv) As the model is in some sense a stochastic version of a spatial system of coupled logistic maps, the restrictions on $m$ in our results are inherited from the behavior of (6): When $m > 4$, the one dimensional deterministic dynamical system (6) would "live" only on a Cantor-like set, and the technique employed in the proof of Theorem 1 would fail. On the other hand, simulations suggest that, even in the case $m > 4$, the random fluctuations can "smooth out" the trajectories so that (4) might survive from initial conditions which would drive (5) to extinction in finitely many steps. The restriction to $m \in (1,3)$ in Theorem 3 stems of course from the fact that this guarantees a unique stable fixed point of the logistic map. It is unclear if Corollary 4 would hold in a situation where (6) has periodic orbits: Then, one can see in simulations large regions of space which are "oscillating out of phase," which might indicate that the system builds up long range structure over large time scales. Similar effects have been studied in [9].

(v) We note that the "stepping stone version of the Bolker–Pacala model" introduced in Definition 1.3 of [10] can be obtained as a scaling limit of a sequence of models considered above: Assume that the parameters of the $N$th model are given by

$$m^{(N)} = 1 + \frac{\alpha M}{N},$$

$$p^{(N)}_{xy} = \frac{1}{N} m_{xy} + \left(1 - \frac{1}{N} \sum_x m_{0x}\right) \delta_{xy},$$

$$\lambda^{(N)}_{xy} = \frac{\alpha \kappa \lambda_{xy}}{N^2},$$

where $\alpha, M, m_{xy}, \lambda_{xy}$ are as in [10], page 191. Let $\xi_0^{(N)}(x) = [N\mu(x)]$, where $\mu$ is some finite measure on $\mathbb{Z}^d$, and define $X_t^{(N)}(x) := \frac{1}{N}\xi_{[Nt]}^{(N)}(x)$. Then $X^{(N)}$ converges in distribution on $D_{[0,\infty)}(\mathcal{M}_f(\mathbb{Z}^d))$ to $X$, the solution of (5) on page 191 of [10], that is, the stepping stone version of the Bolker–Pacala model, with $\gamma = 1$. Of course, this is a remark about finite time horizons. Deducing results about the steady states of the Bolker–Pacala model from our theorems would presumably require a considerable amount of work.

(vi) Hutzenthaler and Wakolbinger [12] have shown that (at least in the case of within-site competition only) the stepping stone version of the Bolker–Pacala model from [10] dies out in any dimension if the carrying capacity, which would correspond to $(m-1)/\sum_x \lambda_{0x}$ in our model, is too small. Similarly, one would expect that our model, even when $m \in (1,3)$, will die out when $\lambda_{xy}$ are too large. Simulations suggest that this is indeed the case, but we have no rigorous proof.

The rest of this paper is organized as follows: in Section 2 we provide a basic lemma showing how "occupancy" spreads through space and prove Theorem 1, in Section 3 we briefly discuss how the results can be transferred to a two-species scenario with (weak) interspecific competition. Section 4 provides results about the deterministic system (5) and proves Theorem 2. These results will be required in Section 5, where we prove Theorem 3 and Corollary 4.

To simplify the notation in the proofs, we will use in the sequel a transformed version of the kernel $\lambda$,

$$\lambda_{xy} = \kappa \gamma_{xy}, \qquad x \neq y, \tag{8}$$

where we assume that $\sum_{y \neq x} \gamma_{xy} = 1$. That is, we separate the nondiagonal part of $\lambda$ into $\kappa := \sum_{x \neq 0} \lambda_{0x}$, the total "nondiagonal" competition and the normalized kernel $\gamma_{0x} = \lambda_{0x}/\kappa$ ($\gamma_{xx} := 0$). Nevertheless, we prefer to state the theorems in terms of $\lambda_{xy}$ because these have an intuitive biological interpretation. For $\eta \in \mathbb{R}_+^{\mathbb{Z}^d}$, $x \in \mathbb{Z}^d$ and $\kappa \geq 0$, we write

$$f_\kappa(x;\eta) := \eta(x)\left(m - \lambda_0 \eta(x) - \kappa \sum_{z \neq x} \gamma_{xz}\eta(z)\right)^+ \tag{9}$$

and

$$F(x;\eta) := \sum_{y \in \mathbb{Z}^d} f_\kappa(y;\eta) p_{yx}. \tag{10}$$

Note that this is just (2) and (3) in the new parametrization.

**2. Survival.** The value $f_\kappa(x;\eta)$ is the mean number of offspring at site $x$ if the present configuration is $\eta$. The maximal (mean) number of offspring at one site in one generation will be denoted by

$$(11) \qquad m^*_{\lambda_0} := \max_{\eta \in \mathbb{R}_+^{\mathbb{Z}^d}} f_\kappa(0;\eta) = \frac{m^2}{4\lambda_0}.$$

If the number of particles at some site $x$ exceeds $M_{\lambda_0} := m/\lambda_0$, then, as the term in the parenthesis in (2) and (9) is negative, no offspring is produced at this site. Furthermore, let us introduce

$$(12) \qquad \bar{m}(\lambda_0, \kappa) := \frac{m-1}{\lambda_0 + \kappa} \quad \text{and} \quad \bar{m}_{\lambda_0} := \bar{m}(\lambda_0, 0),$$

the deterministic equilibrium values when the nondiagonal regulation term is $\kappa$ respectively 0. Note that, for $\eta \equiv \bar{m}(\lambda_0, \kappa)$, we have $f_\kappa(x;\eta) = \bar{m}(\lambda_0, \kappa)$ and, therefore, $\eta(x) = F(x;\eta)$ for all $x \in \mathbb{Z}^d$.

DEFINITION 6. Let $\eta \in \mathbb{R}_+^{\mathbb{Z}^d}$. For a pair of positive numbers $(\varepsilon_1, \varepsilon_2)$, we will say that a site $x$ is $(\varepsilon_1, \varepsilon_2)$-*occupied with respect to* $\eta$ if

$$\eta(x) \in [\varepsilon_1 \bar{m}_{\lambda_0}, (1-\varepsilon_2)M_{\lambda_0}] \quad \text{and} \quad \eta(y) \leq (1-\varepsilon_2)M_{\lambda_0}, \|x-y\|_\infty \leq R_\lambda.$$

We will often say that $\eta(x)$ is $(\varepsilon_1, \varepsilon_2)$-occupied, or just occupied if there is no risk of confusion, meaning that $x$ is $(\varepsilon_1, \varepsilon_2)$-occupied with respect to $\eta$.

As advertised earlier, to prove Theorem 1, we compare the process $(\xi_n)$ with oriented percolation on a sub-grid of $\mathbb{Z}^d \times \mathbb{Z}_+$. The main step is to show that if a site is $(\varepsilon_1, \varepsilon_2)$-occupied with respect to some $\xi_n$, then in a while its neighbors will be also $(\varepsilon_1, \varepsilon_2)$-occupied with high probability. To this end, we consider a perturbed coupled map lattice

$$(13) \qquad \zeta_{n+1}(x) = F(x;\zeta_n) + \delta_n(x),$$

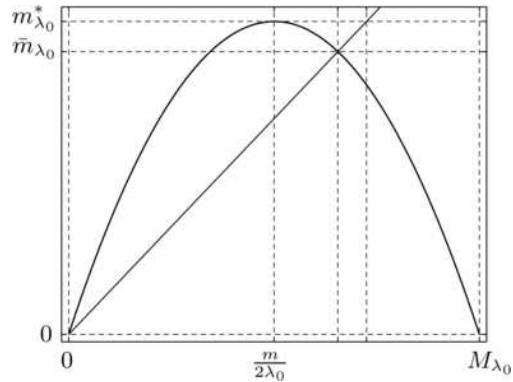

FIG. 2. *Graph of the function* $\tilde{f}(z) = z(m - \lambda_0 z)^+$.

where the perturbation $\delta_n$ is assumed to satisfy $\delta_n(x) \geq -F(x; \zeta_n)$, such that $(\zeta_n)$ is nonnegative. We will show that under certain additional conditions on the perturbation term the system $(\zeta_n)$ has the desired property. Then we view the original process $(\xi_n)$ as a perturbed dynamical system and we will see that the conditions mentioned above are satisfied with high probability if the competition is weak enough.

Let us now introduce and explain some notation which will be used in the sequel. We denote by $p_{xy}^n$ the $n$-step transition probability of a random walk with kernel $p$. As mentioned above, our goal is to show that an occupied site colonizes its neighbors in a couple of steps and remains itself occupied. In the first step the offspring are distributed according to the kernel $p$. Thus, there is in general no reason why an occupied site should remain occupied after one step. Let us fix some $\tilde{m} \in (1, m)$. By the Local Central Limit Theorem, the number

$$(14) \qquad n^* = \min\{j \in \mathbb{N} : p_{0x}^j \tilde{m}^j \geq 1 \text{ for all } x \text{ with } \|x\|_\infty \leq 1\}$$

is finite. We set

$$\mathscr{I} = \{(y, j) \in \mathbb{Z}^d \times \mathbb{Z}_+ : p_{0y}^j > 0, 0 \leq j \leq n^*\}$$
$$\subset \{(y, j) : \|y\|_\infty \leq jR_p, 0 \leq j \leq n^*\}.$$

Suppose that the site 0 is $(\varepsilon_1, \varepsilon_2)$-occupied with respect to $\zeta_0$ and that there is no mass at the other sites. Let us also assume for the moment that the perturbation term vanishes and that the competition between individuals at different sites is zero, that is, $\kappa = 0$. We set $\tilde{f}(z) = z(m - \lambda_0 z)^+$, see Figure 2. As this function is unimodal, to find the minimum of $\tilde{f}$ on some interval, it suffices to consider the values at the endpoints. If for some positive $a$ we have $z \in [a\varepsilon_1 \bar{m}_{\lambda_0}, (1 - \varepsilon_2)M_{\lambda_0}]$, then

$$(15) \quad \tilde{f}(z) \geq \begin{cases} \tilde{f}(\varepsilon_2 M_{\lambda_0}) = \varepsilon_2 M_{\lambda_0} m(1 - \varepsilon_2), & a\varepsilon_1 \bar{m}_{\lambda_0} \geq \varepsilon_2 M_{\lambda_0}, \\ \tilde{f}(a\varepsilon_1 \bar{m}_{\lambda_0}) = a\varepsilon_1 \bar{m}_{\lambda_0} m\left(1 - a\varepsilon_1 + \dfrac{1}{m}\right), & a\varepsilon_1 \bar{m}_{\lambda_0} < \varepsilon_2 M_{\lambda_0}. \end{cases}$$

This means that the number of offspring at site 0 is at least $\tilde{m}\varepsilon_1 \bar{m}_{\lambda_0}$ if $\varepsilon_1$ is sufficiently small. Then the offspring are distributed in the neighborhood according to the kernel $p$. In this neighborhood the mass is again multiplied by at least $\tilde{m}$ (unless the local mass happens to be already large enough) and then distributed according to $p$. Hence, after $k$ steps the mass at a site $y$ is larger than or equal to $p_{0y}^k \tilde{m}^k \varepsilon_1 \bar{m}_{\lambda_0}$. The living space of the whole population at this time is the $k$th timeslice of $\mathscr{I}$ which is contained in the ball of radius $kR_p$. By the definition of $n^*$, after $n^*$ steps the mass in 0 and in points with norm one reaches or maybe exceeds the level $\varepsilon_1 \bar{m}_{\lambda_0}$. Thus, these sites are occupied at that time if the mass there and in the $R_\lambda$-neighborhood does not exceed $(1 - \varepsilon_2)M_{\lambda_0}$.

We need some additional conditions on the perturbation term. We set

(16) $$\mathcal{X} = \{(y,n) \in \mathbb{Z}^d \times \mathbb{Z}_+ : n < n^*, \|y\|_\infty \leq n(R_p + R_\lambda)\}.$$

Consider the assumptions:

$(B1)_{\varepsilon_2}$ for all $(y,n) \in \mathcal{X}$: $F(y;\zeta_n) + \delta_n(y) \leq (1-\varepsilon_2)M_{\lambda_0}$,
$(B2)_{\delta,K}$ for all $(y,n) \in \mathcal{X}$: $F(y;\zeta_n) \geq K$ implies $|\delta_n(y)| \leq \delta F(y;\zeta_n)$.

LEMMA 7. *Assume that $m$ and $p$ are as in Theorem* 1. *For each $K > 0$ and $\delta$ satisfying $m(1-\delta) > \tilde{m} > 1$, there are choices of positive numbers $\varepsilon_1$, $\varepsilon_2$, $\lambda_0^*$ and $\kappa^*$ such that whenever*

(17) $$\lambda_0 \leq \lambda_0^* \quad \text{and} \quad \kappa \leq \kappa^* \lambda_0,$$

*the following holds:*

$\zeta_0(0)$ *is* $(\varepsilon_1, \varepsilon_2)$-*occupied,* $(B1)_{\varepsilon_2}$, $(B2)_{\delta,K}$ *are satisfied*
$$\implies \zeta_{n^*}(x) \text{ are } (\varepsilon_1, \varepsilon_2)\text{-occupied for all } x \text{ with } \|x\|_\infty \leq 1.$$

PROOF. Let $K > 0$ be given. We choose $\varepsilon_2 > 0$ such that

(18) $$m(1-\delta)\left(1 - \varepsilon_2 \frac{m}{m-1}\right) \geq \tilde{m} \quad \text{and} \quad m^*_{\lambda_0} \leq (1-2\varepsilon_2)M_{\lambda_0}.$$

For the second inequality, we need $m < 4$. Then we choose $\varepsilon_1 > 0$ satisfying

(19) $$p_{0y}^n \tilde{m}^n \varepsilon_1 \leq \varepsilon_2 \frac{m+1}{m} \leq \varepsilon_2 \frac{m}{m-1} \quad \text{for all } (n,y) \in \mathcal{I}.$$

Note that this choice guarantees

(20) $$p_{0y}^n \tilde{m}^n \varepsilon_1 \bar{m}_{\lambda_0} \leq \varepsilon_2 \frac{m+1}{m} \bar{m}_{\lambda_0} \leq \varepsilon_2 M_{\lambda_0} \quad \text{for all } (n,y) \in \mathcal{I}.$$

By construction of $\mathcal{I}$, the number $\mathcal{I}_{\min} = \min\{\tilde{m}^n p_{0y}^n : (n,y) \in \mathcal{I}\}$ is positive. Therefore we may choose $\lambda_0^*$ such that, for $\lambda_0 \leq \lambda_0^*$,

$$\varepsilon_1 \bar{m}_{\lambda_0} \mathcal{I}_{\min} \geq K.$$

Finally, we choose $\kappa^*$ such that, for some $\alpha$ satisfying $\tilde{m} < m - \alpha$,

$$(1-\varepsilon_2)M_{\lambda_0} \kappa^* \lambda_0^* \leq \alpha.$$

Let us first consider the case $\kappa = 0$. We have to show that

$$\zeta_{n^*}(x) \in [\varepsilon_1 \bar{m}_{\lambda_0}, (1-\varepsilon_2)M_{\lambda_0}], \quad \|x\|_\infty \leq 1.$$

By $(B1)_{\varepsilon_2}$, we have

(21) $$\zeta_{n+1}(x) = F(x;\zeta_n) + \delta_n(x) \leq (1-\varepsilon_2)M_{\lambda_0} \quad \text{for all } (x,n) \in \mathcal{X}.$$

This means, in particular, $\zeta_{n^*}(x) \leq (1-\varepsilon_2)M_{\lambda_0}$ for $\|x\|_\infty \leq 1$. To complete the proof for that case, we show by induction on $n$ that

$$\zeta_n(y) \in [p_{0y}^n \tilde{m}^n \varepsilon_1 \bar{m}_{\lambda_0}, (1-\varepsilon_2)M_{\lambda_0}], \qquad 0 \leq n \leq n^* \quad \text{and} \quad (y,n) \in \mathscr{I}.$$

By definition of $n^*$, the assertion of the lemma then follows. For $n=0$, the claim holds by assumption. If it holds for some $n < n^*$, then, first using (20) and (15), then (19) and the first part of (18), we obtain

$$(1-\delta)f(y;\zeta_n) \geq (1-\delta)\tilde{f}(p_{0y}^n \tilde{m}^n \varepsilon_1 \bar{m}_{\lambda_0})$$
$$\geq (1-\delta)p_{0y}^n \tilde{m}^n \varepsilon_1 \bar{m}_{\lambda_0} \cdot m\left(1 - p_{0y}^n \tilde{m}^n \varepsilon_1 + \frac{1}{m}\right)$$
$$\geq p_{0y}^n \tilde{m}^{n+1} \varepsilon_1 \bar{m}_{\lambda_0}.$$

Hence,

$$(1-\delta)F(y;\zeta_n) = \sum_{z \in \mathbb{Z}^d}(1-\delta)f(z;\zeta_n)p_{zy}$$
$$\geq \sum_{z \in \mathbb{Z}^d} p_{0z}^n \tilde{m}^{n+1} \varepsilon_1 \bar{m}_{\lambda_0} p_{zy} = \varepsilon_1 \bar{m}_{\lambda_0} \tilde{m}^{n+1} p_{0y}^{n+1}, \qquad (y,n) \in \mathscr{I}.$$

In particular, we have $F(y;\zeta_n) \geq \varepsilon_1 \bar{m}_{\lambda_0} \tilde{m}^{n+1} p_{0y}^{n+1} \geq K$ for $\lambda_0 \leq \lambda_0^*$. Therefore, $(B2)_{\delta,K}$ applies and from the last display we obtain

$$\zeta_{n+1}(y) \geq (1-\delta)F(y;\zeta_n) \geq \varepsilon_1 \bar{m}_{\lambda_0} \tilde{m}^{n+1} p_{0y}^{n+1}.$$

This concludes the proof of the induction and proves the lemma in the special case $\kappa = 0$.

Now let us turn to the case $\kappa > 0$. Assumption $(B1)_{\varepsilon_2}$, (17) and the choice of $\kappa^*$ imply that

$$0 \leq \kappa \sum_{y \neq x} \gamma_{xy}\zeta_n(y) \leq \kappa^* \lambda_0 (1-\varepsilon_2)M_{\lambda_0} \leq \alpha,$$
$$\|x\|_\infty \leq n(R_\lambda + R_p) - R_\lambda, n < n^*,$$

where $\alpha > 0$ satisfies $m - \alpha > \tilde{m}$. We obtain

$$f^{(l)}(x;\zeta_n) := \zeta_n(x)(m - \alpha - \lambda_0\zeta_n(x))^+ \leq f_\kappa(x;\zeta_n(x))$$
$$\leq \zeta_n(x)(m - \lambda_0\zeta_n(x))^+ =: f^{(u)}(x;\zeta_n).$$

So we can use the same induction as in the diagonal case. For the lower bound estimates, we use $f^{(l)}$ and for the upper bound estimates, we use $f^{(u)}$. □

We set $\zeta_0 = \xi_0$ and assume that $(\zeta_n)$ is the solution of (13) with the perturbation term

$$\delta_n(x) = N^{(n,x)}(F(x;\xi_n)) - F(x;\xi_n).$$

Thus, $(\xi_n)$ with $\xi_n = \zeta_n$ can be considered as a perturbed coupled map lattice.

PROOF OF THEOREM 1. Recall the definition of the space-time box in (16). For $(x,n) \in \mathbb{Z}^d \times \mathbb{Z}_+$, we set

$$X(x,n) = \{N^{(y,j)} : (y,j) \in (x,n) + \mathcal{X}\}.$$

Consider the events

$$A(x,n) = \{N^{(y,j)}(m^*_{\lambda_0}) \leq (1-\varepsilon_2)M_{\lambda_0}, (y,j) \in (x,n) + \mathcal{X}\}$$

and

$$B(x,n) = \left\{ \sup_{(y,j) \in (x,n) + \mathcal{X}} \sup_{t \geq K} \left| \frac{N^{(y,j)}(t)}{t} - 1 \right| \leq \delta \right\}.$$

We say that $X(x,n)$ is *good* if $A(x,n) \cap B(x,n)$ holds. First we want to show that the probability of a good realization can be made arbitrarily large by choosing small $\lambda_0$. It is of course enough to consider the corresponding problem in the space-time point $(0,0)$. As $A(0,0)$ implies $(B1)_{\varepsilon_2}$ and $B(0,0)$ implies $(B2)_{\delta,K}$ on the event $A(0,0) \cap B(0,0)$, Lemma 7 yields

$$\{\xi_0(0)(\varepsilon_1,\varepsilon_2)\text{-occupied}\} \cap (A(0,0) \cap B(0,0))$$
$$\subset \{\xi_{n^*}(y), \|y\|_\infty \leq 1 \ (\varepsilon_1,\varepsilon_2)\text{-occupied}\}.$$

By translation invariance, the corresponding statement is also true for all $(x,n) \in \mathbb{Z}^d \times \mathbb{Z}_+$. Furthermore, we point out that $X(x,n)$ and $X(x',n')$ are independent if $\|x-x'\|_\infty \geq 2n(R_\lambda + R_p)$ or $|n-n'| > n^*$.

Let $\Delta$ be the number of points in $\mathcal{X}$ and let $(N(t))_{t \geq 0}$ be a standard Poisson process. Then we have

$$\mathbb{P}[A(0,0)] = (1-a(\lambda_0))^\Delta \qquad \text{where } a(\lambda_0) = \mathbb{P}[N(m^*_{\lambda_0}) > (1-\varepsilon_2)M_{\lambda_0}].$$

According to (18), we have $m^*_{\lambda_0} \leq (1-2\varepsilon_2)M_{\lambda_0}$. Thus, for some $\tilde{c}_1 > 0$, we have

$$a(\lambda_0) = \mathbb{P}\left[\frac{N(m^*_{\lambda_0})}{m^*_{\lambda_0}} - 1 > \frac{(1-\varepsilon_2)M_{\lambda_0}}{m^*_{\lambda_0}} - 1\right]$$
$$\leq \mathbb{P}\left[\frac{N(m^*_{\lambda_0})}{m^*_{\lambda_0}} - 1 > \varepsilon_2\right] \leq \exp\left(-\frac{\tilde{c}_1 \varepsilon_2^2}{\lambda_0}\right).$$

Furthermore, by standard large deviation results for Poisson processes, for some $\tilde{c}_2 > 0$ and sufficiently large $K$, we have

$$\mathbb{P}[B(0,0)] = \mathbb{P}\left[\sup_{t \geq K}\left|\frac{N(t)}{t} - 1\right| \leq \delta\right]^\Delta = \left(1 - \mathbb{P}\left[\sup_{t \geq K}\left|\frac{N(t)}{t} - 1\right| > \delta\right]\right)^\Delta$$
$$\geq (1 - \exp(-\tilde{c}_2 \delta^2 K))^\Delta.$$

From the proof of Lemma 7, one can see that making $K$ large corresponds to making $\lambda_0$ small. Hence,

$$(22) \qquad \mathbb{P}[(A(0,0) \cap B(0,0))^c] \leq \mathbb{P}[A(0,0)^c] + \mathbb{P}[B(0,0)^c] \leq \theta(\lambda_0),$$

where $\theta(\lambda_0) \leq \exp(-c/\lambda_0)$ for some suitable positive constant $c = c(p, m, R_\lambda)$. This implies

$$\mathbb{P}[X(0,0) \text{ is good}] \geq 1 - \theta(\lambda_0) = 1 - (1 - \sqrt{p(\lambda_0)})^\Delta,$$

where $p(\lambda_0) = (1 - \theta(\lambda_0)^{1/\Delta})^2$. Since $p(\lambda_0)$ converges to one as $\lambda_0$ goes to 0 and the range of dependence is finite, it is clear that the good sites percolate if $\lambda_0$ is small enough. For example, one can apply a result by Liggett, Schonmann and Stacey (see [16], Theorem 26) to show that, for fixed $n$, the distribution of the random field $\mathbb{1}_{\{X(x,n) \text{ is good}\}}$ dominates the product measure $\nu_{p(\lambda_0)} = \bigotimes_{\mathbb{Z}^d} \text{Ber}(p(\lambda_0))$ on $\{0,1\}^{\mathbb{Z}^d \times \mathbb{Z}_+}$. Comparison of the process $(\mathbb{1}_{\{X(x,n) \text{ is good}\}})_{x \in \mathbb{Z}^d \times n^* \mathbb{Z}_+}$ with independent oriented percolation concludes the proof. $\square$

**3. A competing species model.** In this section we consider two processes $(\xi_n^{(1)})$ and $(\xi_n^{(2)})$, modeling, for example, two different species or genetic types living in the same habitat and competing for similar (or the same) resources. In the absence of the other type, each of them is a version of the basic process described in the introduction, possibly with different parameters.

Let $(\lambda_{xy}^{(ij)})_{x,y \in \mathbb{Z}^d}$, $i,j \in \{1,2\}$, be translation invariant nonnegative kernels on $\mathbb{Z}^d$ with finite range $R_\lambda$. These kernels will determine the intra- respectively interspecific competition: The average reproductive success of an $i$-individual at $x$ is reduced by each $j$-individual at $y$ by $\lambda_{xy}^{(ij)}$. The evolution of $(\xi_n^{(1)}, \xi_n^{(2)})$ may then be described as follows. Similar to the single species model, we define

$$f_1(x; \xi_n^{(1)}, \xi_n^{(2)}) = \xi_n^{(1)}(x) \left( m_1 - \sum_y \lambda_{xy}^{(11)} \xi_n^{(1)}(y) - \sum_y \lambda_{xy}^{(12)} \xi_n^{(2)}(y) \right)^+,$$

$$f_2(x; \xi_n^{(1)}, \xi_n^{(2)}) = \xi_n^{(2)}(x) \left( m_2 - \sum_y \lambda_{xy}^{(22)} \xi_n^{(2)}(y) - \sum_y \lambda_{xy}^{(21)} \xi_n^{(1)}(y) \right)^+,$$

$$F_1(x; \xi_n^{(1)}, \xi_n^{(2)}) = \sum_y f_1(y; \xi_n^{(1)}, \xi_n^{(2)}) p_{yx}^{(1)},$$

$$F_2(x; \xi_n^{(1)}, \xi_n^{(2)}) = \sum_y f_2(y; \xi_n^{(1)}, \xi_n^{(2)}) p_{yx}^{(2)},$$

where $m_i$ is the mean number of offspring of a type $i$ individual in the absence of competition. If $N_1^{(x,n)}$, $N_2^{(x,n)}$, $(x,n) \in \mathbb{Z}^d \times \mathbb{Z}_+$ are independent standard Poisson processes on $\mathbb{R}_+$, then, given $(\xi_n^{(1)}, \xi_n^{(2)})$, the configuration of the next generation is given by

$$(\xi_{n+1}^{(1)}, \xi_{n+1}^{(2)}) = (N_1^{(x,n)}(F_1(x; \xi_n^{(1)}, \xi_n^{(2)})), N_2^{(x,n)}(F_2(x; \xi_n^{(1)}, \xi_n^{(2)}))).$$

We obtain the following about long-term coexistence if the competition terms are weak enough:

THEOREM 8. *For given $m_1, m_2 \in (1,4)$, $p^{(1)}$ and $p^{(2)}$ satisfying* (A1) *and range $R_\lambda$, there are positive numbers $\lambda_1^*$, $\lambda_2^*$, $\kappa_1^*$, $\kappa_2^*$ and $\gamma^*$ such that if the conditions:*

(i) $0 < \lambda_0^{(ii)} \leq \lambda_0^*$, $\sum_{y \neq x} \lambda_{xy}^{(ii)} \leq \lambda_0^{(ii)} \kappa_i^*$, $i \in \{1, 2\}$;
(ii) $\sum_y \lambda_{xy}^{(12)}, \sum_y \lambda_{xy}^{(21)} \leq \gamma^* \min\{\lambda_0^{(11)}, \lambda_0^{(22)}\}$;

*are satisfied, then both populations survive with positive probability, provided that for some $x, y \in \mathbb{Z}$ we have $f_1(x; \xi_0^{(1)}, \xi^{(2)}(0)) > 0$ and $f_2(y; \xi_0^{(11)}, \xi_0^{(22)}) > 0$. Furthermore, conditioned on survival of both populations,*

$$\liminf_{N \to \infty} \frac{1}{N} \sum_{n=1}^{N} \mathbb{1}_{\{\xi_n^{(1)}(0)\xi_n^{(2)}(0) > 0\}} > 0 \qquad a.s.,$$

*that is, we have local coexistence.*

To prove this theorem, one can essentially use the same argument as we have used in the proof of Lemma 7 to reduce the case $\kappa > 0$ to the case $\kappa = 0$.

**4. Results for the deterministic system.** In this section we will prove Theorem 2. For clarity of exposition, we start with the "diagonal case" $\kappa = 0$. Let us consider more generally a coupled map lattice $(\zeta_n)$ on $\mathbb{Z}^d$, defined via

(23) $$\zeta_{n+1}(x) = \sum_{y \in \mathbb{Z}^d} g(\zeta_n(y)) p_{yx}, \qquad x \in \mathbb{Z}^d,$$

where $(p_{yx})_{x,y \in \mathbb{Z}^d}$ is a translation invariant stochastic kernel with finite range satisfying (A1) and $g : [0, G] \to [0, G]$ is a continuously differentiable function. We think of the single site function $g$ as having 0 as a repelling fixed point and another stable fixed point $\bar{a} \in (0, G]$ which attracts $(0, G]$, that is, for any $x_0 \in (0, G]$, the sequence $(x_n)$ defined through $x_{n+1} = g(x_n)$ converges to $\bar{a}$. [Thus, in particular, $g'(0) > 1$, $g(G) > 0$.] Then obviously $\zeta \equiv 0$ and $\zeta \equiv \bar{a}$ are fixed points of (23), and one is strongly inclined to believe that in this well-behaved scenario there are no others. We will say that a dynamical system $(\eta_n)$ on $\mathbb{Z}^d$ *converges locally* to $a \in \mathbb{R}$ if for each finite $\Lambda \subset \mathbb{Z}^d$ and each $\varepsilon > 0$ there exists $N_0$ such that

$$|\eta_n(x) - a| \leq \varepsilon \qquad \text{for all } x \in \Lambda \text{ and } n \geq N_0.$$

Having been unable to find the result we need in the literature, we provide Lemma 9 below. Assume the following:

(DS1) For each $a > 0$, there exist sequences $(\alpha_n)$ and $(\beta_n)$ such that $0 < \alpha_0 \leq a$, $\beta_0 = G$, $\alpha_n \uparrow \bar{a}$, $\beta_n \downarrow \bar{a}$ and $g([\alpha_n, \beta_n]) \subset [\alpha_{n+1}, \beta_{n+1}]$.

Note that this implies the following:

(DS2) There exists $a \in (0, \bar{a})$ with the following property:
If $\zeta_0(0) \in [a, G]$, then there is $N_0 \in \mathbb{N}$ such that $\zeta_{N_0}(x) \in [a, G]$, $\|x\|_\infty \leq 1$.

A proof that (DS1) $\Rightarrow$ (DS2) is basically a reformulation of the proof of Lemma 7. Note that (DS1) holds true, for example, if we assume additionally that $g$ is concave (see, e.g., the construction given in Lemma 12). We refrain from pursuing the most general conditions for (DS1), but observe that this together with Lemma 9 already yields a proof of Theorem 2 in the diagonal case $\kappa = 0$.

LEMMA 9. *If $\zeta_0(x) \in (0, G]$ for some $x \in \mathbb{Z}^d$ and (DS1) holds, then $(\zeta_n)$ converges locally to $\bar{a}$.*

In the following we will call the set $\mathcal{N}_k(A) := \{x \in \mathbb{Z}^d : \inf_{y \in A} \|x - y\|_\infty \leq k\}$ the $k$-neighborhood of $A$. If $A = \{x\}$, then we write $\mathcal{N}_k(x)$ for the $k$-neighborhood of $x$.

PROOF OF LEMMA 9. Let $\Lambda$ be a finite subset of $\mathbb{Z}^d$. We may assume that $\Lambda$ is a ball with respect to the sup norm. Let $(\alpha_n)$ and $(\beta_n)$ be sequences from (DS1). Given $\varepsilon > 0$, we choose $n_0$ such that $\beta_n - \alpha_n < \varepsilon$ holds for all $n \geq n_0$. According to (DS2), there exist $a \in (0, \bar{a})$ and $n_1 \in \mathbb{N}$ such that

$$\zeta_n(x) \geq a \Rightarrow \zeta_{n+n_1}(y) \geq a \qquad \text{for all } y \text{ with } \|x - y\|_\infty \leq 1.$$

Since 0 and $\bar{a}$ are the only fixed points of $g$, $g'(0) > 1$ and $a \in (0, \bar{a}]$, we have $g(a) \geq a$. It follows that if for all $y$ in the $R_p$-neighborhood of some point $x$ we have $\zeta_0(y) \geq a$, then

$$\zeta_1(x) = \sum_y g(\zeta_0(y)) p_{yx} \geq a.$$

We set

$$\Lambda' := \mathcal{N}_{R_p(n_0+n_1)}(\Lambda) \quad \text{and} \quad \Lambda_i = \mathcal{N}_{R_p(n_0-i)}(\Lambda), \qquad i \in \{0, \ldots, n_0\}.$$

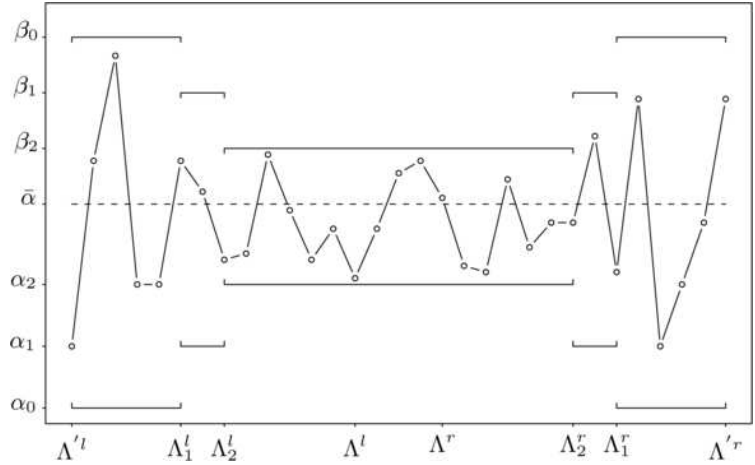

Fig. 3. *A sketch of the nested boxes and the dynamical system in one dimension after two steps. The superscripts $l$ and $r$ indicate the left respectively the right-hand side of the box. Initially the population size is in $[\alpha_0, \beta_0]$ in each site of $\Lambda'$. After two steps the population size is in $[\alpha_0, \beta_0]$ on $\Lambda' \setminus \Lambda_1$, is in $[\alpha_1, \beta_1]$ on $\Lambda_1 \setminus \Lambda_2$ and in $[\alpha_2, \beta_2]$ on $\Lambda_2$.*

Note that $\Lambda_i = \mathcal{N}_{R_p}(\Lambda_{i+1})$ and $\Lambda_{n_0} = \Lambda$. By (DS2), there is some time point $n_2 \in \mathbb{N}$ such that $\zeta_{n_2}(x) \geq a$ for all $x \in \Lambda'$. We claim that $\zeta_{n_2+n}(x) \geq a$ for all $x \in \Lambda_0$ and all $n \geq 0$. Indeed, during the next $n_1 - 1$ steps from time $n_2$ on, the mass in the points of the $R_p$-neighborhood of $\Lambda_0$ remain bounded from below by $a$. According to (DS2) by the time $n_2 + n_1$ each point in the 1-neighborhood of $\Lambda'$ is bounded below by $a$. Hence, we are, in particular, again in the above situation.

For simplicity of notation, we assume that $\zeta_0(x) \geq a$ for all $x \in \Lambda'$. We need to show that $\zeta_n(x) \in [\alpha_{n_0}, \beta_{n_0}]$ for all $x \in \Lambda$ and $n \geq n_0$. First, we check inductively that for $n = 0, 1, \ldots, n_0$ we have the following (see Figure 3 for an illustration):

$$
\begin{aligned}
&\text{(i)} && \zeta_n(x) \in [\alpha_n, \beta_n] && \text{for } x \in \Lambda_n, \\
&\text{(ii)} && \zeta_k(x) \in [\alpha_k, \beta_k] && \text{for } x \in \Lambda_k \setminus \Lambda_{k+1}, k = 0, 1, \ldots, n-1.
\end{aligned}
\tag{24}
$$

For $n = 0$, (i) is true by assumption, and (ii) is void. Assume that (i) and (ii) hold true for some $n < n_0$, let $k \in \{0, 1, \ldots, n+1\}$ and $x \in \Lambda_k \setminus \Lambda_{k-1}$, respectively $x \in \Lambda_{n+1}$ if $k = n+1$. As $\Lambda_{k-1} = \mathcal{N}_{R_p}(\Lambda_k)$, we have

$$\zeta_{n+1}(x) = \sum_y g(\zeta_n(y)) p_{yx} = \sum_{y \in \Lambda_{k-1}} g(\zeta_n(y)) p_{yx} \in [\alpha_k, \beta_k]$$

by (DS1), proving (i) and (ii) for $n + 1$.

To conclude the proof, note that, by the argument above, the set of configurations $\zeta$ such that

$$\zeta(x) \geq a \qquad \text{for } x \in \Lambda',$$
$$\zeta(x) \in [\alpha_k, \beta_k] \qquad \text{for } x \in \Lambda_k \setminus \Lambda_{k+1}, k = 0, 1, \ldots, n_0 - 1,$$
$$\zeta(x) \in [\alpha_{n_0}, \beta_{n_0}] \qquad \text{for } x \in \Lambda_{n_0}$$

is invariant under the dynamics (23), hence, we have, in particular, for $n \geq n_0$,

$$\zeta_n(x) \in [\alpha_{n_0}, \beta_{n_0}] \qquad \text{for } x \in \Lambda (= \Lambda_{n_0}). \qquad \square$$

For the "nondiagonal" case $\kappa > 0$, we need three more lemmas. Note that we only need to consider the case $\lambda_0 = 1$. Otherwise, consider $\tilde{\zeta}$ defined by $\tilde{\zeta}_n(x) = \lambda_0 \zeta_n(x)$, which solves the iteration given by (9) and (10) with $\lambda_0$ replaced by 1 and $\kappa$ by $\kappa/\lambda_0$. Until the end of this section we write $\bar{m}_{1,\kappa} = \bar{m}(1, \kappa)$, $m^* = m_1^* = m^2/4$ and $\bar{m} = \bar{m}_1 = m - 1$ [see (12)].

LEMMA 10. *There exist positive $\kappa^*$ and $\delta$ such that for $\kappa \leq \kappa^*$ exist sequences $(\alpha_n)$, $(\beta_n)$ in $[\bar{m}_{1,0} - \delta, \bar{m}_{1,0} + \delta]$ satisfying the following:*

1. $\alpha_n \uparrow \bar{m}_{1,\kappa}$, $\beta_n \downarrow \bar{m}_{1,\kappa}$;
2. *If $\zeta(y) \in [\alpha_n, \beta_n]$ for all $y \in \mathcal{N}_{R_\lambda}(x)$, then $f_\kappa(x; \zeta) \in [\alpha_{n+1}, \beta_{n+1}]$.*

PROOF. For fixed $x \in \mathbb{Z}^d$, we may consider the mapping $\zeta \mapsto f_\kappa(x; \zeta)$ as a function of the restriction of $\zeta$ to the $R_\lambda$-neighborhood of $x$ [viewed as an element of $\mathbb{R}^k$ where $k$ is the number of points in $\mathcal{N}_{R_\lambda}(x)$]. We denote by $\vec{m}_{1,\kappa}$ the vector of length $k$ with all entries equal to $\bar{m}_{1,\kappa}$ and by $B_\delta(\vec{m}_{1,\kappa})$ the $\delta$-neighborhood of $\vec{m}_{1,\kappa}$ with respect to sup norm.

The gradient of $\zeta \mapsto f_\kappa(x; \zeta)$ is given by [we assume that the positive part appearing in (9) is not 0]

$$\partial_{\zeta(x)} f_\kappa(x; \zeta) = m - 2\zeta(x) - \kappa \sum_{y \neq x} \gamma_{xy} \zeta(y),$$

$$\partial_{\zeta(y)} f_\kappa(x; \zeta) = -\kappa \gamma_{xy} \zeta(x) \qquad \text{for } y \neq x.$$

Choose positive $\varepsilon$, $\delta$ and $\kappa^*$ satisfying

$$(|m - 2| + 2\delta + \kappa^*(\delta + m - 1))^2 < 1 - \varepsilon \quad \text{and}$$
(25)
$$\kappa^* < \min\left\{\frac{\delta}{m-1}, \frac{\sqrt{\varepsilon}}{\sqrt{2}(2+\delta)}\right\}.$$

For $\zeta \in B_\delta(\vec{m}_{1,0})$, we have

$$\partial_{\zeta(x)} f_\kappa(x; \zeta) \leq m - 2(\bar{m}_{1,0} - \delta) - \kappa \sum_{y \neq x} \gamma_{xy}(\bar{m}_{1,0} - \delta)$$
$$= 2 - m + 2\delta - \kappa(m-1) + \kappa\delta$$

and

$$\partial_{\zeta(x)} f_\kappa(x;\zeta) \geq m - 2(\bar{m}_{1,0} + \delta) - \kappa \sum_{y \neq x} \gamma_{xy}(\bar{m}_{1,0} + \delta)$$
$$= 2 - m - 2\delta - \kappa(m-1) - \kappa\delta,$$

hence,

$$|\partial_{\zeta(x)} f_\kappa(x;\zeta)| \leq |m - 2| + (m-1)\kappa + 2\delta + \kappa\delta$$

and due to (25), we obtain, for $\kappa \leq \kappa^*$

(26) $$(\partial_{\zeta(x)} f_\kappa(x;\zeta))^2 < 1 - \varepsilon.$$

For $y \neq x$, we have

$$|\partial_{\zeta(y)} f_\kappa(y;\zeta)| = \kappa \gamma_{xy}\zeta(x) \leq \kappa \gamma_{xy}(\bar{m}_0 + \delta)$$
$$\leq \kappa\gamma_{xy}(m-1) + \delta\gamma_{xy}\kappa < \kappa\gamma_{xy}(2+\delta).$$

Consequently,

(27) $$\sum_{y \neq x}(\partial_{\zeta(y)} f_\kappa(y;\zeta))^2 < (2+\delta)^2 \kappa^2 \sum_{x \neq y} \gamma_{xy}^2 \leq (2+\delta)^2 \kappa^2 < \frac{\varepsilon}{2},$$

where the last inequality holds if (25) is satisfied.

Altogether, the above implies that for all $\zeta \in B_\delta(\vec{m}_{1,\kappa})$ and $\kappa \leq \kappa^*$ we have

(28) $$\|\nabla f_\kappa(x;\zeta)\|_2^2 < 1 - \frac{\varepsilon}{2}.$$

Due to the mean value theorem for all $\zeta, \zeta' \in B_\delta(\vec{m}_{1,0})$ exists $\tilde{\zeta} \in B_\delta(\vec{m}_{1,0})$ such that

$$|f_\kappa(x;\zeta) - f_\kappa(x;\zeta')| = |\nabla f_\kappa(x;\tilde{\zeta})(\zeta - \zeta')|$$
$$\leq \|\nabla f_\kappa(x;\tilde{\zeta})\|_2 \cdot \|\zeta - \zeta'\|_2 \leq c\|\zeta - \zeta'\|_2,$$

where $c = \sqrt{1 - \varepsilon/2} < 1$. Thus, the claim of the lemma follows. We only need to note that $f_\kappa(x;\vec{m}_{1,\kappa}) = m_{1,\kappa}$ and that $|\bar{m}_{1,0} - \bar{m}_{1,\kappa}| < \delta$ if $\kappa < \delta/(m-1)$ which holds by (25). □

LEMMA 11. *For each $\delta > 0$ exists $\kappa^* > 0$ such that whenever $\kappa \leq \kappa^*$ and $f_\kappa(x;\zeta_0) > 0$ for some $x \in \mathbb{Z}^d$, the following holds: For each finite $\Lambda \subset \mathbb{Z}^d$, there exists $N \in \mathbb{N}$ such that $\zeta_n(x) \in [\bar{m}_{1,0} - \delta, \bar{m}_{1,0} + \delta]$ for all $x \in \Lambda$ and all $n \geq N$.*

PROOF. Recall our assumption $\lambda_0 = 1$, which implies $M_{\lambda_0} = m$. For all $x \in \mathbb{Z}^d$, $\zeta \in [0,m]^{\mathbb{Z}^d}$ and $\tilde{\delta} > 0$, we have

$$\kappa \sum_{z \neq x} \gamma_{xz} \zeta(z) \leq m\kappa.$$

That implies

$$f_{\kappa,l}(\zeta(x)) := \zeta(x)(m - m\kappa - \zeta(x)) \leq f_\kappa(x;\zeta) \leq \zeta(x)(m - \zeta(x)) =: f_u(\zeta(x)).$$

The nonzero fixed points of $f_{\kappa,l}$ and $f_u$ are respectively $\bar{m}_l = m - m\kappa - 1$ and $\bar{m}_{1,0}$. Furthermore, if $m\kappa < \delta$, then $\bar{m}_{1,0} - \bar{m}_l < \delta$.

According to Lemma 7, there is $n_1 \in \mathbb{N}$ and $a > 0$ with the property

(29) $$\zeta_n(x) \geq a \Rightarrow \zeta_{n+n_1}(y) \geq a, \qquad \|x - y\|_\infty \leq 1.$$

Thus, for each finite $\Lambda' \subset \mathbb{Z}^d$, there is $n_2 \in \mathbb{N}$ such that $\zeta_{n_2}(x) \geq a$ for all $x \in \Lambda'$.

According to Lemma 12 for each $\delta > 0$, one can choose $\kappa^*$ and sequences $(a_n)$ and $(b_n)$ such that for all $\kappa \leq \kappa^*$ the following holds:

$$a_0 \leq a,$$
$$f_{\kappa,l}([a_n, b_n]), f_u([a_n, b_n]) \subset [a_{n+1}, b_{n+1}],$$
for some $n_0 \in \mathbb{N}$: $\quad a_n, b_n \in [\bar{m}_0 - \delta, \bar{m}_0 + \delta]$ for all $n \geq n_0$.

A construction analogous to the proof of Lemma 9 concludes the proof. $\square$

The following lemma is a deterministic ingredient in our construction [see (DS1)], providing a shrinking sequence of intervals which the one-point iteration maps into themselves. Having been unable to find a proof in the literature, we provide one here. The property in question will hold for a concave $f$ with 0 a repelling and another attracting fixed point and does not depend on the particular functional form of $f$. On the other hand, as we also need to consider a slightly perturbed version $f_\delta$ (where in our case the perturbation is of a particular functional form), we refrain from generality and stick to $f_\delta, f : [0,m] \to [0,m^*]$,

(30) $$f_\delta(x) = x(m - \delta - x)^+, \qquad f(x) = x(m - x)^+,$$

where $m^* = m^2/4 = \max f = f(m/2)$. Recall that $\bar{m} = m - 1$, $\bar{m}_\delta = m - \delta - 1$ are the (unique) attracting fixed points of $f$ respectively $f_\delta$ (we think of small $\delta$).

LEMMA 12. *Letting $m \in (1,3)$, consider $f, f_\delta$ as defined in* (30). *For each $\varepsilon > 0$, one can choose positive $\gamma$ and $\tilde{\varepsilon}$, a strictly increasing sequence $(\alpha_n)$, and a strictly decreasing sequence $(\beta_n)$ with the following properties:*

(A) *There exists* $N_0 \in \mathbb{N}$ *s.t.* $\alpha_n, \beta_n \in [\bar{m} - \varepsilon, \bar{m} + \varepsilon]$ *for all* $n \geq N_0$.
(B) *For all* $n \leq N_0$ *and* $0 \leq \delta \leq \gamma$: $f_\delta([\alpha_n, \beta_n]), f([\alpha_n, \beta_n]) \subset [\frac{\alpha_{n+1}}{1-\tilde{\varepsilon}}, \frac{\beta_{n+1}}{1+\tilde{\varepsilon}}]$.

*Furthermore, $\alpha_0 > 0$ can be chosen arbitrarily small and $\beta_0 < m$ can be chosen arbitrarily close to $m$.*

PROOF. We wish to construct the sequences $(\alpha_n)$ and $(\beta_n)$ in such a way that

$$\alpha_n < \alpha_{n+1} < \bar{m}_\gamma \leq \bar{m} < \beta_{n+1} < \beta_n \tag{31}$$

and

$$f_\gamma([\alpha_n, \beta_n]), f([\alpha_n, \beta_n]) \subset (\alpha_{n+1}, \beta_{n+1}) \tag{32}$$

for all $n$. This together with

$$\bar{m} - \varepsilon < \lim_{n \to \infty} \alpha_n \leq \lim_{n \to \infty} \beta_n < \bar{m} + \varepsilon \tag{33}$$

will suffice to conclude, as $f_\gamma(x) \leq f_\delta(x) \leq f(x)$ for $0 \leq \delta \leq \gamma$ and (32) implies (B) for each finite $N_0$ and sufficiently small $\tilde{\varepsilon}$. The construction is slightly different depending on whether the slope of $f$ at its attractive fixed point $\bar{m}$ is $\in (0, 1)$, $= 0$ or $\in (-1, 0)$, thus, we consider the cases $m \in (1, 2)$, $m = 2$ and $m \in (2, 3)$ separately.

Letting $m \in (1, 2)$, choose $\gamma \in (0, \varepsilon)$ s.t. $m - \gamma \in (1, 2)$. Take arbitrary $\alpha_0 \in (0, \bar{m} - \gamma)$ and $\beta_0 > m/2$ s.t. $f_\gamma(\beta_0) \geq f_\gamma(\alpha_0)$. This guarantees $f([\alpha_0, \beta_0])$, $f_\gamma([\alpha_0, \beta_0]) \subset [f_\gamma(\alpha_0), m^*]$. Define

$$\alpha_{n+1} = \frac{\alpha_n + f_\gamma(\alpha_n)}{2}, \qquad n \geq 0,$$

$$\beta_1 = \frac{m^* + \frac{m}{2}}{2} \quad \text{and} \quad \beta_{n+1} = \frac{f(\beta_n) + \beta_n}{2}, \qquad n \geq 1.$$

Note that $m^* < m/2$ in the case considered, so the choice of $\beta_1$ ensures (32) for $n = 0$ and that $f, f_\gamma$ are increasing on $[\alpha_0, \beta_1]$. As $f_\gamma(x) > x$ on $(0, \bar{m}_\gamma)$ and $f'_\gamma(\bar{m}_\gamma) \geq 0$, we have $\alpha_n < \alpha_{n+1} < f_\gamma(\alpha_n)$ for $n \geq 1$. Thus, $\alpha_n \nearrow \bar{m}_\gamma$. Similarly, observing that $x > f(x) \geq \bar{m}$ for $x \in (\bar{m}, m/2)$, we have $\beta_n > \beta_{n+1} > f_\gamma(\beta_n)$ for $n \geq 1$, hence, $\beta_n \searrow \bar{m}$. This proves (31), (32) and (33) in this case.

Let $m = 2$. In this case $f(m/2) = m^*$, so the values of $f(\beta_n)$ cannot be decreasing, and we modify the construction as follows: Choose $0 < \gamma < \varepsilon$. Picking $\alpha_0 \in (0, \gamma)$, define

$$\alpha_{n+1} = \frac{f_\delta(\alpha_n) + \alpha_n}{2},$$

$$\beta_n = \frac{2-\gamma}{2} + \sqrt{(2-\gamma)^2/4 - f_\gamma(\alpha_n)}, \qquad n = 0, 1, \ldots$$

As above, we have $\alpha_n \nearrow \bar{m}_\gamma = 2 - \gamma$. Note that $\beta_n$ is the larger root of $f_\gamma(x) = f_\gamma(\alpha_n)$, and that the solutions of $f_\gamma(x) = \bar{m}_\gamma$ are $\bar{m}_\gamma = 1 - \gamma$ and 1 in the case $m = 2$, so that $f(\beta_n) \geq f_\gamma(\beta_n) > \alpha_{n+1}$ and $\beta_n \searrow 1$. Hence, (31), (32) and (33) are satisfied.

Finally, let $m \in (2, 3)$. Here, as $\bar{m} > m/2$, we need to observe that $f([\alpha_n, \beta_n])$ will contain $m^*$ as long as $\alpha_n \leq m/2$, so $\beta_n$ must not decrease too quickly. Furthermore, as $f'(\bar{m}) < 0$, once we come close to $\bar{m}$, the roles of the lower and upper boundary are interchanged in each step.

Choose $\gamma > 0$ s.t. $m - \gamma \in (2, 3)$ and $\bar{m}_\gamma > f_\gamma(m^*) > m/2$. Pick $\alpha_0 \in (0, (m - \gamma)/2)$. While $(\alpha_n + f_\gamma(\alpha_n))/2 \leq m/2$, we set

$$\alpha_{n+1} = \frac{\alpha_n + f_\gamma(\alpha_n)}{2}.$$

Let $n_0$ be the smallest integer satisfying $(\alpha_{n_0} + f_\gamma(\alpha_{n_0}))/2 > m/2$. We set

$$\alpha_{n_0+1} = \frac{\alpha_{n_0} + f_\gamma(\alpha_{n_0})}{2} \wedge \frac{1}{2}\left(\frac{m}{2} + f_\gamma(m^*)\right).$$

Now we choose $\beta_0, \ldots, \beta_{n_0}$ s.t. $m^* < \beta_i < \beta_{i-1} < m$ and $f_\gamma(\beta_i) > \alpha_{i+1}$, $i = 1, \ldots, n_0$. Note that this is possible because $f_\gamma(m^*) > m/2$. Put $\beta_{n_0+1} = (\beta_{n_0} + m^*)/2$.

Let us check (31) and (32) for $n \leq n_0$: as $f_\gamma(x) > x$ for $x \in (0, \bar{m}_\gamma)$ and $f_\gamma(m^*) < \bar{m}_\gamma$, the sequence $(\alpha_n)_{n \in \{0, \ldots, n_0+1\}}$ is strictly increasing. $(\beta_n)_{n \in \{0, \ldots, n_0+1\}}$ is strictly decreasing by construction. By definition, we have

$$f_\gamma(\alpha_n) \geq 2\alpha_{n+1} - \alpha_n > \alpha_{n+1}.$$

Note that while $\alpha_n \leq m/2$, that is, $n \leq n_0$, we always have

$$f_\gamma([\alpha_n, \beta_n]), f([\alpha_n, \beta_n]) \subset (\alpha_{n+1}, m^*] \subset (\alpha_{n+1}, \beta_{n+1}).$$

For $n \geq n_0 + 1$, define

(34) $\qquad \alpha_{n+1} = \frac{1}{2}(f_\gamma(\beta_n) + \alpha_n), \qquad \beta_{n+1} = \frac{1}{2}(\beta_n + f(\alpha_n)).$

In order to verify (31) and (32) for $n \geq n_0 + 1$, consider

(35) $\qquad a \in \left(\frac{m}{2}, \bar{m}_\gamma\right), b \in (\bar{m}, m) \qquad$ satisfying $f(a) < b, f_\gamma(b) > a$.

Note that then

$$a' = \tfrac{1}{2}(a + f_\gamma(b)) \quad \text{and} \quad b' = \tfrac{1}{2}(b + f(a))$$

fulfill

$$a' \in (a, \bar{m}_\gamma), b' \in (\bar{m}, b) \quad \text{and} \quad f(a') < b', f_\gamma(b') > a'.$$

Indeed, by assumption, we have $f_\gamma(b) > a$, so $a' > a$. On the other hand, $f_\gamma(b) < \bar{m}_\gamma$ because $f_\gamma$ is decreasing in $[\bar{m}_\gamma, m]$ and $b > \bar{m}_\gamma = f_\gamma(\bar{m}_\gamma)$. As $f$ is decreasing in the considered region, we have

$$f(a') < f(a) < \tfrac{1}{2}(b + f(a)) = b'.$$

Similarly, $b' \in (\bar{m}, b)$ and $f_\gamma(b') > a'$.

Obviously, $a = \alpha_{n_0+1}$ and $b = \beta_{n_0+1}$ satisfy the condition (35), hence, (31) and (32) hold true for $n > n_0$ as well.

By the above construction, $\alpha_n \nearrow \alpha \in (m/2, \bar{m}_\gamma]$, $\beta_n \searrow \beta \in [\bar{m}, m^*)$, where $(\alpha, \beta)$ solves $f(\alpha) = \beta$, $f_\gamma(\beta) = \alpha$. For $\gamma = 0$, the unique solution would be $\alpha = \beta = \bar{m}$, for $\gamma$ sufficiently small, we have (33). $\square$

PROOF OF THEOREM 2. Let $\Lambda$ be a finite ball in $\mathbb{Z}^d$ and $\varepsilon > 0$. Let $n_1 \in \mathbb{N}$ be such that (29) is fulfilled. For the sequences $(\alpha_n)$, $(\beta_n)$ from Lemma 10, choose $n_0$ s.t. $\beta_n - \alpha_n \leq \varepsilon$ for all $n \geq n_0$. Define $\Lambda'$ and $\Lambda_0, \ldots, \Lambda_{n_0}$ through

$$\Lambda' = \mathcal{N}_{(n_0+n_1)(R_\lambda + R_p)}(\Lambda) \quad \text{and} \quad \Lambda_i = \mathcal{N}_{(n_0-i)(R_\lambda + R_p)}(\Lambda), i \in \{0, \ldots, n_0\}.$$

According to Lemma 11, there exists $n_2 \in \mathbb{N}$ such that $\zeta_n(x) \in [\bar{m}_{1,0} - \delta, \bar{m}_{1,0} + \delta]$ for all $x \in \Lambda_0$ and $n \geq n_2$. Then, for simplicity of notation, we may assume $n_2 = 0$. Now the rest of the proof is a reproduction of the arguments from the proof of Lemma 9. $\square$

**5. Coupling.** In this section we prove Theorem 3 and Corollary 4. Let us first describe the idea behind the successful coupling. Recall in (7) the definition of the coupling $(\xi^{(1)}, \xi^{(2)})$. Consider three large (but finite) boxes $B_1 \subset B_2 \subset B_3 \subset \mathbb{Z}^d$ and assume that $\xi^{(1)}$ and $\xi^{(2)}$ agree on $B_1$ with values close to $\bar{m}_{\lambda_0}$, that they are close to $\bar{m}_{\lambda_0}$ but do not necessarily agree on $B_2$, and that on $B_3$ all sites are occupied in both systems. In view of Lemma 7, we expect that the region of sites which are occupied in both systems grows. If the competition is not too strong, the random system "follows closely" the deterministic one. Thus, in view of Theorem 2, we can hope that the region where both systems are close to the deterministic equilibrium $\bar{m}_{\lambda_0}$ is growing as well. Finally, there is a chance that Poisson variables whose means are close to each other produce the same realization. Therefore, there is also hope that the region where both systems are the same grows too.

Thus, for suitably tuned parameters, we expect that, with high probability, the above situation will reproduce itself after some time on larger boxes $B_1' \subset B_2' \subset B_3'$. As before, this observation lends itself to a comparison with finite range dependent percolation on a coarse grained space-time grid. A certain subtlety stems from the problem that the coarse graining must be chosen depending on $\lambda_0$ in such a way that the dependence range of the percolation does not diverge when taking $\lambda_0$ small.

For $k, l \in \mathbb{N}$, we set $A_k = \mathcal{N}_{k(R_\lambda+R_p+1)}(0)$ and $A_{k,l} = \mathcal{N}_{k(R_\lambda+R_p+1)+l}(0)$. Letting $\mathcal{X}(y, n)$, $(y, n) \in \mathbb{Z}^d \times \mathbb{Z}_+$ be the event that for some $N \in \mathbb{N}$, to be chosen later, the following holds:

$$\xi_n^{(1)}(x) = \xi_n^{(2)}(x) \in \left[\frac{m-1-\delta}{\lambda_0}, \frac{m-1+\delta}{\lambda_0}\right] =: I(m, \delta, \lambda_0)$$
$$\text{for all } x \in y + A_N,$$

(36) $\quad \xi_n^{(1)}(x), \xi_n^{(2)}(x) \in I(m, \delta, \lambda_0) \qquad \text{for all } x \in y + A_{4N} \setminus A_N,$

$$\xi_n^{(1)}(x), \xi_n^{(2)}(x) \in [\varepsilon_1 \bar{m}_{\lambda_0}, (1-\varepsilon_2) M_{\lambda_0}] =: J(m, \lambda_0)$$
$$\text{for all } x \in y + A_{7N} \setminus A_{4N}.$$

Our goal is to show that the process $\mathbb{1}_{\mathcal{X}(y,n)}$ dominates oriented independent percolation on a suitable sub-grid of $\mathbb{Z}^d \times \mathbb{Z}_+$. The main part of the proof is carried out in Lemma 13 below. With this lemma one can, for example, use the Liggett, Schonmann and Stacey argument as we have done in the proof of Theorem 1.

Let $n^*$ be as defined in (14) and note that this number only depends on $m$ and on the kernel $p$. As we will later choose $N$ large, we will be able to choose it as a multiple of $n^*$. In the sequel we will assume that $N/n^*$ is an integer.

LEMMA 13. *For $m \in (1,3)$, $p$ as in assumption* (A1) *and $\tilde{\varepsilon} > 0$, there exist $\lambda_0^*, \kappa^* > 0$ such that for each $\lambda_0 \leq \lambda_0^*$, $\kappa \leq \kappa^* \lambda_0$ one can choose $N$ such that*

(37) $\quad \mathbb{P}[\mathcal{X}(y, n+N) \text{ for all } y \text{ with } \|x-y\|_\infty \leq N/n^* | \mathcal{X}(x,n)] \geq 1 - \tilde{\varepsilon}$

*holds for all $x \in \mathbb{Z}^d$.*

PROOF. Let $m \in (1,3)$ and $\tilde{\varepsilon} > 0$ be given. Due to translation invariance and the Markov property, the left-hand side in (37) does not depend on $(x, n)$. Thus, it is enough to prove

(38) $\quad \mathbb{P}[\mathcal{X}(y, N) \text{ for all } y \text{ with } \|y\|_\infty \leq N/n^* | \mathcal{X}(0,0)] \geq 1 - \tilde{\varepsilon}.$

Choose positive $\varepsilon$, $\delta$ and $\kappa^*$ satisfying

(39) $\quad |m - 2| + 2\delta + \kappa^*(\delta + m - 1) < 1 - \varepsilon \quad \text{and}$

$$\kappa^* < \min\left\{\frac{\delta}{m-1}, \frac{\varepsilon}{2(2+\delta)}\right\}.$$

These constants also satisfy (25). Thus, the properties of $f_\kappa$ [see (9)], proven in Lemma 10, are preserved. Note that, unlike the situation in Lemma 10, we do not set $\lambda_0 = 1$ here. Furthermore, similar to (26), (27) and (28), we obtain

$$\|\nabla f_\kappa(x;\zeta)\|_1 \leq 1 - \frac{\varepsilon}{2} \tag{40}$$

$$\text{if } \zeta(y) \in \left[\frac{m-1-\delta}{\lambda_0}, \frac{m-1+\delta}{\lambda_0}\right] \quad \text{for all } y \in \mathcal{N}_{R_\lambda}(x).$$

We choose $k_0$ such that, for all $k \geq k_0$, we have

$$\frac{|A_{k+1}|}{|A_k|}\left(1 - \frac{\varepsilon}{2}\right) \leq \frac{|A_{k_0+1}|}{|A_{k_0}|}\left(1 - \frac{\varepsilon}{2}\right) =: c(\varepsilon) < 1. \tag{41}$$

We will assume that $N \geq k_0$. We set $\mathcal{X}_0 = \mathcal{X}(0,0)$ and $\mathcal{X}_N = \mathcal{X}_{N,1} \cap \mathcal{X}_{N,2} \cap \mathcal{X}_{N,3}$, where

$$\mathcal{X}_{N,1} = \{\xi_N^{(1)}(x) = \xi_N^{(2)}(x) \in I(m,\delta,\lambda_0) \text{ for all } x \in A_{3N}\},$$

$$\mathcal{X}_{N,2} = \{\xi_N^{(1)}(x), \xi_N^{(2)}(x) \in I(m,\delta,\lambda_0) \text{ for all } x \in A_{6N} \setminus A_{3N}\},$$

$$\mathcal{X}_{N,3} = \{\xi_N^{(1)}(x), \xi_N^{(2)}(x) \in J(m,\lambda_0) \text{ for all } x \in A_{7N,N/n^*} \setminus A_{6N}\}.$$

Furthermore, we define for each $n \leq N$ the event $\Psi_n$ by

$$\Psi_n = \left\{\forall (x,k) \in \bigcup_{j=1}^n A_{4(N-j)} \times \{j\} : \xi_k^{(1)}(x), \xi_k^{(2)}(x) \in I(m,\delta,\lambda_0)\right\}.$$

As $\mathcal{X}_N$ implies that $\mathcal{X}(y,N)$ holds for all $y$ with $\|y\|_\infty \leq N/n^*$, $\mathbb{P}[\mathcal{X}_N|\mathcal{X}_0]$ is a lower bound for the left-hand side of (38). Therefore, it suffices to show $\mathbb{P}[\mathcal{X}_N^c|\mathcal{X}_0] \leq \tilde{\varepsilon}$. Because

$$\mathbb{P}[\mathcal{X}_N^c|\mathcal{X}_0] \leq \mathbb{P}[\mathcal{X}_{N,1}^c \cap \Psi_N|\mathcal{X}_0] + \mathbb{P}[\mathcal{X}_{N,2}^c \cap \Psi_N|\mathcal{X}_0] \tag{42}$$
$$+ \mathbb{P}[\mathcal{X}_{N,3}^c|\mathcal{X}_0] + \mathbb{P}[\Psi_N^c|\mathcal{X}_0],$$

it suffices to estimate each of the summands. To do this, we will repeatedly use large deviation estimates for Poisson random variables. There are constants $c_1$ and $\delta_1$ such that

$$\mathbb{P}[\Psi_N^c|\mathcal{X}_0] \leq N|A_{4N}|\exp\left(-\frac{c_1\delta_1^2}{\lambda_0}\right). \tag{43}$$

Now let us consider the first term on the right-hand side of (42). We denote $\mathcal{F}_m = \sigma(\{N_j^{(x,l)} : x \in \mathbb{Z}^d, l \leq m\}, \xi_0^{(1)}, \xi_0^{(2)})$. We have

$$\frac{1}{|A_{3N}|} \sum_{x \in A_{3N}} \mathbb{E}[|\xi_N^{(1)}(x) - \xi_N^{(2)}(x)|\mathbb{1}_{\Psi_N}|\mathcal{F}_{N-1}]$$

$$\leq \mathbb{1}_{\Psi_{N-1}} \frac{1}{|A_{3N}|} \sum_{x \in A_{3N}} \mathbb{E}[|\xi_N^{(1)}(x) - \xi_N^{(2)}(x)||\mathcal{F}_{N-1}]$$

$$\leq \mathbb{1}_{\Psi_{N-1}} \frac{1}{|A_{3N}|} \sum_{x \in A_{3N}} \sum_{y \in \mathcal{N}_{R_p}(x)} p_{yx} \sum_{z \in \mathcal{N}_{R_\lambda}(y)} |\nabla_z f_\kappa(y;\tilde{\xi})||\xi^{(1)}_{N-1}(z) - \xi^{(2)}_{N-1}(z)|$$

$$\leq \mathbb{1}_{\Psi_{N-1}} \sum_{z \in A_{3N+1}} |\xi^{(1)}_{N-1}(z) - \xi^{(2)}_{N-1}(z)| \frac{1}{|A_{3N}|} \sum_{y \in \mathcal{N}_{R_\lambda}(z)} |\nabla_z f_\kappa(y;\tilde{\xi})| \sum_{x \in A_{3N}} p_{xy}$$

$$\leq \mathbb{1}_{\Psi_{N-1}} \frac{|A_{3N+1}|}{|A_{3N}|} \left(1 - \frac{\varepsilon}{2}\right) \frac{1}{|A_{3N+1}|} \sum_{z \in A_{3N+1}} |\xi^{(1)}_{N-1}(z) - \xi^{(2)}_{N-1}(z)|$$

$$\leq \mathbb{1}_{\Psi_{N-1}} c(\varepsilon) \frac{1}{|A_{3N+1}|} \sum_{z \in A_{3N+1}} |\xi^{(1)}_{N-1}(z) - \xi^{(2)}_{N-1}(z)|.$$

We can iterate the above argument to obtain on $\mathcal{X}_0$

(44)
$$\frac{1}{|A_{3N}|} \sum_{x \in A_{3N}} \mathbb{E}[|\xi^{(1)}_N(x) - \xi^{(2)}_N(x)| \mathbb{1}_{\Psi_N} | \mathcal{F}_0]$$
$$\leq c(\varepsilon)^N \frac{1}{|A_{4N}|} \sum_{z \in A_{4N}} |\xi^{(1)}_0(z) - \xi^{(2)}_0(z)|$$
$$\leq c(\varepsilon)^N \frac{1}{|A_{4N}|} \sum_{z \in A_{4N} \setminus A_N} 2\delta \bar{m}_{\lambda_0}$$
$$= c(\varepsilon)^N \frac{|A_{4N} \setminus A_N|}{|A_{4N}|} 2\delta \bar{m}_{\lambda_0} \leq c(\varepsilon)^N 2\delta \bar{m}_{\lambda_0}.$$

From this, we obtain

(45)
$$\mathbb{P}[\mathcal{X}^c_{1,N} \cap \Psi_N | \mathcal{X}_0] \leq \sum_{x \in A_{3N}} \mathbb{E}[|\xi^{(1)}_N(x) - \xi^{(2)}_N(x)| \mathbb{1}_{\Psi_N}]$$
$$\leq c(\varepsilon)^N 2\bar{m}_{\lambda_0} \delta |A_{3N}|.$$

Note that on $\mathcal{X}_0$ for all $|x| \leq R_\lambda + R_p$ we have $\xi^{(1)}_n(x) = \xi^{(2)}_n(x)$ for all $n \leq N-1$.

To estimate the second term of the right-hand side of (42), let $(\alpha_n)$ and $(\beta_n)$ be sequences from Lemma 12 satisfying $\alpha_0 \leq \varepsilon_1(m-1)$ and $\beta_0 \geq (1-\varepsilon_2)m$. Let $\kappa^*$ be small enough for Theorem 1 and Theorem 2 to apply. Let $N_0$ be the number from Lemma 12 such that, for all $n \geq N_0$, we have $\alpha_n/((1-\tilde{\delta})\lambda_0), \beta_n/((1+\tilde{\delta})\lambda_0) \in I(m,\lambda_0,\delta)$. Recall that in the formulation of Lemma 12 we have chosen $\lambda_0 = 1$, but it holds for general $\lambda_0$. We assume $N_0 \leq N$. If for all $x \in \mathcal{N}_{R_\lambda+R_p}(0)$ we have $\xi_0(x) \in [\alpha_n/\lambda_0, \beta_n/\lambda_0]$, where $\xi$ is a version of the processes considered, then there exist positive constants $c_2$ and $\delta_2$ such that, for all $n \leq N_0$, we have

$$\mathbb{P}\left[\xi_1(0) \notin \left[\frac{\alpha_{n+1}}{\lambda_0}, \frac{\beta_{n+1}}{\lambda_0}\right]\right] = \mathbb{P}\left[N^{(0,0)}(F(0;\xi_0)) \notin \left[\frac{\alpha_{n+1}}{\lambda_0}, \frac{\beta_{n+1}}{\lambda_0}\right]\right]$$

$$\leq \exp\left(-\frac{c_2\delta_2^2}{\lambda_0}\right),$$

because $F(0;\xi_0) \in [\alpha_{n+1}/((1-\tilde{\delta})\lambda_0), \beta_{n+1}/((1+\tilde{\delta})\lambda_0)]$. It follows that

$$(46) \qquad \mathbb{P}[\mathcal{X}_{N,2}^c \cap \Psi_N | \mathcal{X}_0] \leq N|A_{7N} \setminus A_{4N}| \exp\left(-\frac{c_2\delta_2^2}{\lambda_0}\right).$$

The upper bound for the third term on the right-hand side of (42) is obtained as follows:

$$\begin{aligned}
\mathbb{P}[\mathcal{X}_{N,3}^c | \mathcal{X}_0] &= \mathbb{P}[\exists x \in A_{7N,N/n^*} : \xi_N^{(1)}(x), \xi_N^{(2)}(x) \notin J(m,\lambda_0) | \mathcal{X}_0] \\
&\leq \mathbb{P}\bigg[\exists k \in \left\{1,\ldots,\frac{N}{n^*}-1\right\} \\
&\qquad \exists x \in A_{7N,k} \setminus A_{6N} : \xi_{kn^*}^{(1)}(x) \text{ or } \xi_{kn^*}^{(2)}(x) \notin J(m,\lambda_0) | \mathcal{X}_0\bigg] \\
&\leq \frac{2N|A_{7N,N/n^*-1} \setminus A_{6N}|}{n^*}\theta(\lambda_0),
\end{aligned}$$

(47)

where $\theta(\lambda_0) \leq \exp(-c_3/\lambda_0)$ for some positive $c_3$ is defined in (22).

Let $N$ be the smallest multiple of $n^*$ larger than $1/\lambda_0$. Using the above estimates, one can choose some positive $c$ and $r \in \mathbb{N}$ such that

$$\mathbb{P}[\mathcal{X}_N^c | \mathcal{X}_0] \leq \exp\left(-\frac{c}{\lambda_0}\right) N^r.$$

The right-hand side goes to zero as $\lambda_0$ goes to zero. Thus, (38) follows. □

Before we turn to the proof of Theorem 3, we need a result about oriented percolation. Let $\theta \in (0,1)$ be given and let $A(x,n)$, $(x,n) \in \mathbb{Z}^d \times \mathbb{Z}_+$ be i.i.d. Bernoulli random variables with parameter $\theta$. For $k < n$, we say that $(x,k)$ *is connected to* $(y,n)$, this will be denoted by $(x,k) \to (y,n)$, if there is a sequence $x = x_0, \ldots, x_{n-k} = y$ such that $\|x_i - x_{i-1}\|_\infty \leq 1$ and $A(x_i, k+i) = 1$ for $i = 1, \ldots, n-k$. Let $C_0 = \{(x,n) : (0,0) \to (x,n)\}$ be the cluster of the origin. We call a space time-point $(y,n)$ $C_0$-*exposed* if there exists a sequence $y_n, \ldots, y_0$ such that $y_n = y$, $\|y_k - y_{k-1}\|_\infty \leq 1$, and $(y_k, k) \notin C_0$, $k = 1, \ldots, n$.

The next lemma follows from [6]. The idea behind the proof is as follows: With the "usual" percolation interpretation in mind, let us call a site $(x,n)$ *wet* if there is a backward path $(x,n) = (x_0,n), (x_1,n-1), \ldots, (x_n,0)$ with $\|x_i - x_{i-1}\| \leq 1$ consisting only of open sites, that is, $A(x_i, n-i) = 1$, $i = 0, 1, \ldots, n-1$. Otherwise, the site will be called *dry*. Lemma 7 in [6] shows, using a contour-counting argument, that if $\theta$ is sufficiently close to 1, the dry sites do not percolate. In fact, this lemma even obtains an exponential bound on the tail of the size of the cluster of dry sites containing a given site. The next ingredient is complete convergence for oriented percolation ([6],

Lemma 8): When $\theta$ is close enough to 1, there is a fixed $c > 0$ and a random $N_0$ such that on $\{|C_0| = \infty\}$, $\{(x,n) : (x,n) \text{ wet and } \|x\| \leq cn\} \subset C_0$ for all $n \geq N_0$. In words, any wet site inside the "cone" $\{(x,n) : \|x\| \leq cn, n \geq N_0\}$ is also connected to $(0,0)$ by an open path. Fix $c' \in (0,c)$. Assume that $\{|C_0| = \infty\}$, consider $(y,n)$ with $\|y\| \leq c'n$ and $n \geq 2N_0$, say. If $(y,n)$ is $C_0$-exposed, there must be a backward path $(y,n) = (y_0, n), (y_1, n-1), \ldots, (y_n, 0)$ with $(y_i, n-i) \notin C_0$. By the above, at least the initial $n(c-c')/2$ of these sites must be dry [for otherwise, they would be in $C_0$, as they must satisfy $\|y_i\| \leq c(n-i)$]. Hence, there must be a cluster of dry sites containing a point in $\{(x,n) : \|x\| \leq c'n\}$ of size at least $n(c-c')/2$. By the exponential bound on the cluster size distribution and the Borel–Cantelli lemma, this does not occur for $n$ sufficiently large.

LEMMA 14. *If $\theta$ is sufficiently close to 1, then there is a positive constant $c$ such that, for large enough times $n$ conditional on $\{|C_0| = \infty\}$, there are no $C_0$-exposed sites in $\{x \in \mathbb{Z}^d : \|x\|_\infty \leq cn\}$.*

PROOF OF THEOREM 3.  Recall the definition of the event $\mathcal{X}(y,n)$ from (36). Theorem 1 implies that, conditioned on nonextinction of $(\xi_n^{(1)})$ and $(\xi_n^{(2)})$, with probability one, there exist some finite time $N_0$ such that the event $\mathcal{X}(0, N_0)$ holds. Therefore, we may assume a priori that $\mathcal{X}(0,0)$ holds.

We set $\widetilde{N} = [N/(2n^*)]$, $B = \{(x,n) \in \mathbb{Z}^d \times \mathbb{Z}_+ \mid \|x\|_\infty \leq N, n \leq N\}$, $L = \widetilde{N}\mathbb{Z}^d$ and $K = N\mathbb{Z}_+$. Then we have

$$\mathbb{Z}^d \times \mathbb{Z}_+ = \bigcup_{(\alpha,\nu) \in L \times K} ((\alpha,\nu) + B).$$

Let $\|\cdot\|_L$ be the norm on $L$ defined by $\|\alpha\|_L = \|\alpha\|_\infty / \widetilde{N}$. To prove the theorem, it is enough to show that for each $x^* \in \mathbb{Z}^d$ there is time $T$, such that $\xi_n^{(1)}(x^*) = \xi_n^{(2)}(x^*)$ holds for all $n \geq T$. Let us fix an arbitrary $x^* \in \mathbb{Z}^d$ and let $\alpha^* \in L$ be such that $\|\alpha^* - x^*\|_\infty \leq \widetilde{N}$. We define a process $(\eta_\nu)$ on the coarse-grained lattice $L \times K$ by

$$\eta_0(\alpha) = \mathbb{1}_{\mathcal{X}(\alpha,0)} \quad \text{and} \quad \eta_\nu(x) = \mathbb{1}_{\mathcal{X}(\alpha, \nu-N)}, \qquad \nu > 0.$$

Note that $\mathbb{1}_{\mathcal{X}(\alpha,\nu-N)} = 1$ for $\nu > 0$ ensures that $\xi_k^{(1)}(y) = \xi_k^{(2)}(y)$ holds for all $(y,k) \in (\alpha, \nu - N) + B$, because any backward in time path starting in $(y,k)$ will at time $\nu - N$ be inside $\alpha + A_N$, where $\xi^{(1)}$ and $\xi^{(2)}$ are the same on the event $\mathcal{X}(\alpha, \nu - N)$. In particular, $\eta_\nu(\alpha^*) = 1$ implies $\xi_k^{(1)}(x^*) = \xi_k^{(2)}(x^*)$ for all $k \in \{\nu - N, \ldots, \nu\}$. We aim at showing that, for suitable choice of parameters, the process $(\eta_\nu)$ dominates oriented percolation on $L \times K$. To this end, we need to estimate

$$\mathbb{P}[\eta_{\nu+N}(\beta) = 1, \|\alpha - \beta\|_L \leq 1 | \eta_\nu(\alpha) = 1],$$

whereas, due to translation invariance, it is enough to consider the corresponding probability for $(\alpha, \nu) = (0, 0)$. By the construction of $(\eta_\nu)$ and Lemma 13 for each positive $\tilde{\varepsilon}$, one can choose $\lambda_0$, $\kappa$ and $N$ such that

$$\mathbb{P}[\eta_N(\beta) = 1, \|\beta\|_L \leq 1 | \eta_0(0) = 1]$$
$$\geq \mathbb{P}\bigg[\mathcal{X}(z, N), \|z\|_\infty \leq \frac{N}{n^*} \bigg| \mathcal{X}(0, 0)\bigg] \geq 1 - \tilde{\varepsilon}.$$

From the proofs of Theorem 1 and Lemma 13, it can be seen that for $x$ with $\|x\|_\infty \leq 1$ the event $\mathcal{X}(\widetilde{N}x, N)$ is independent of the Poisson processes [which generate $(\xi_n^{(1)})$ and $(\xi_n^{(2)})$] outside the box $\{(y, k) \in \mathbb{Z}^d \times \mathbb{Z}_+ : k \leq N, \|y\|_\infty \leq (8N+2)(R_\lambda + R_p)\}$. Therefore, $(\eta_\nu)$ can be considered as $M$-dependent oriented percolation on $L \times K$, where $M = 20n^*(R_\lambda + R_p) \geq (8N + 2)(R_\lambda + R_p)/\widetilde{N}$. Note that $M$ does not depend on $N$ and $\lambda_0$. Thus, the fact that we need to make $\lambda_0$ small does not affect the comparison.

Let $\theta$ be close enough to 1 such that Lemma 14 holds. For $\tilde{\varepsilon} \in (0, (1 - \sqrt{\theta})^\Delta)$, where $\Delta = |\{(\alpha, \nu) \in L \times K : \nu \in \{0, N\}, \|\alpha\|_L \leq M\}|$, we have

$$\mathbb{P}[\eta_N(\beta) = 1, \|\beta\|_L \leq 1 | \eta_0(0) = 1] \geq 1 - (1 - \sqrt{\theta})^\Delta.$$

As in the proof of Theorem 1, according to Theorem B26 in [15], $(\eta_\nu)$ dominates the nearest neighbor oriented percolation build from the product measure $\nu_\theta$ on $\{0, 1\}^{L \times K}$. Thus, we obtain $\mathbb{P}[|C_\eta| = \infty] > 0$, where $C_\eta \subset L \times K$ is the cluster of the origin generated by $(\eta_\nu)$. By Lemma 14, conditioned on $\{|C_\eta| = \infty\}$, there is a time $T$ such that the points $(\alpha^*, \nu) \in L \times K$ with $\nu \geq T$ are not $C_\eta$-exposed.

We claim that, for each $n \geq T$, conditioned on $\{|C_\eta| = \infty\}$, we have $\xi_n^{(1)}(x^*) = \xi_n^{(2)}(x^*)$. If we assume the contrary, then there must be a path $(x^*, n) = (x_n, n), (x_{n-1}, n-1), \ldots, (x_0, 0)$ in $\mathbb{Z}^d \times \mathbb{Z}_+$ such that $\|x_{i+1} - x_i\|_\infty \leq R_\lambda + R_p$ and $\xi_i^{(1)}(x_i) \neq \xi_i^{(2)}(x_i)$ for all $i \in \{0, \ldots, n-1\}$. From this path, we discard the points $(x_i, i)$ for which $i$ is not a multiple of $N$, thus obtaining for some integer $k$ the path $(x_{kN}, kN), (x_{(k-1)N}, (k-1)N), \ldots, (x_0, 0)$. To this path belongs a path $(\alpha^*, (k+1)N), (\alpha_{kN}, kN), \ldots, (\alpha_0, 0)$ in $L \times K$ where for $j \in \{1, \ldots, k\}$ we choose $\alpha_{jN}$ such that $(x_{(j-1)N}, (j-1)N) \in (\alpha_{jN}, jN) + B$ and $\alpha_0$ such that $\|\alpha_0 - x_0\|_\infty \leq N$. The assumption means that $\eta_{iN}(\alpha_{iN}) = 0$ for all $i \in \{0, \ldots, k\}$. This contradicts the fact that $(\alpha^*, (k+1)N)$ is not $C_\eta$-exposed. $\square$

PROOF OF COROLLARY 4. The sequence $(\xi_n)$, seen as a sequence of random measures on $\mathbb{Z}^d$, is relatively compact with respect to convergence in distribution in the vague topology because the expectation of $\xi_n(x)$ is bounded uniformly by $m_{\lambda_0}^*$.

It is clear that Dirac measure in $\mathbf{0} \in \mathbb{Z}_+^{\mathbb{Z}^d}$ is invariant. If there were two invariant distributions assigning probability 0 to the configuration $\mathbf{0}$, then

Theorem 3 would imply that they coincide on finite subsets of $\mathbb{Z}^d$ and, therefore, they must be equal.

It remains to prove the existence of a limiting invariant distribution $\mu$ satisfying $\mu(\mathbf{0}) = 0$. Let the initial distribution $\mu_0$ be the product measure on $\mathbb{Z}^d$ such that $\xi_0(x) = N^{(0,x)}(\bar{m}(\lambda_0, \kappa))$ for all $x \in \mathbb{Z}^d$. Let $\mu_n$ be the distribution of $\xi_n$. Then the Cesaro average $1/N \sum_{n=0}^{N} \mu_n$ converges along some subsequence $\{N_k\}$ to some measure $\bar{\mu}$. This measure is invariant for $(\xi_n)$ (see, e.g., [14], Proposition I.1.8).

To show $\bar{\mu}(\mathbf{0}) = 0$, it is enough to prove that the restriction of $(\xi_n)$ to $\mathbb{Z}$ survives with probability 1. At time 0, each site is occupied in the sense of Definition 6 with probability

$$\mathbb{P}[N^{(0,0)}(\bar{m}(\lambda_0, \kappa)) \in [\varepsilon_1 \bar{m}_{\lambda_0}, (1-\varepsilon_2)M_{\lambda_0}]],$$

where $\varepsilon_1$ and $\varepsilon_2$ are as in the proof of Lemma 7. In particular, at time 0, there are infinitely many occupied sites. Again, by comparison with oriented percolation, we have $\mathbb{P}_{\xi_0}[\xi_n = \mathbf{0} \text{ for some } n] = 0$ because supercritical percolation starting from infinitely many wet sites does not die out (see, e.g., Theorem B24 in [15]). $\square$

## REFERENCES

[1] BLATH, J., ETHERIDGE, A. M. and MEREDITH, M. E. (2007). Coexistence in locally regulated competing populations and survival of branching annihilating random walk. *Ann. Appl. Probab.* To appear.
[2] BOLKER, B. M. and PACALA, S. W. (1999). Spatial moment equations for plant competition: Understanding spatial strategies and the advantages of short dispersal. *American Naturalist* **153** 575–602.
[3] BOLKER, B. M., PACALA, S. W. and NEUHAUSER, C. (2003). Spatial dynamics in model plant communities: What do we really know? *American Naturalist* **162** 135–148.
[4] CHAZOTTES, J.-R. and FERNANDEZ, B., eds. (2005). *Dynamics of Coupled Map Lattices and of Related Spatially Extended Systems. Lecture Notes in Physics* **671**. Springer, New York.
[5] DOBRUSHIN, R. L., KRYUKOV, V. I. and TOOM, A. L., eds. (1990). *Stochastic Cellular Systems*: *Ergodicity, Memory, Morphogenesis*. Manchester Univ. Press.
[6] DURRETT, R. (1992). Multicolor particle systems with large threshold and range. *J. Theoret. Probab.* **5** 127–152. MR1144730
[7] DURRETT, R. (1999). Stochastic spatial models. *SIAM Rev.* **41** 677–718. MR1722998
[8] DURRETT, R. and LEVIN, S. A. (1994). Stochastic spatial models: A user's guide to ecological applications. *Philos. Trans. Roy. Soc. London. Ser. B* **343** 329–350. MR1678307
[9] DURRETT, R. and LEVIN, S. (2000). Lessons on pattern formation from planet WATOR. *J. Theor. Biol.* **205** 201–214.
[10] ETHERIDGE, A. M. (2004). Survival and extinction in a locally regulated population. *Ann. Appl. Probab.* **14** 188–214. MR2023020
[11] FOURNIER, N. and MÉLÉARD, S. (2004). A microscopic probabilistic description of a locally regulated population and macroscopic approximations. *Ann. Appl. Probab.* **14** 1880–1919. MR2099656

[12] HUTZENTHALER, M. and WAKOLBINGER, A. (2007). Ergodic behavior of locally regulated branching populations. *Ann. Appl. Probab.* **17** 474–501.
[13] LAW, R. and DIECKMANN, U. (2002). Moment approximations of individual-based models. In *The Geometry of Ecological Interactions* (U. Dieckmann, R. Law and J. A. J. Metz, eds.) 252–270. Cambridge Univ. Press.
[14] LIGGETT, T. M. (1985). *Interacting Particle Systems*. Springer, New York. MR0776231
[15] LIGGETT, T. M. (1999). *Stochastic Interacting Systems*: *Contact, Voter and Exclusion Processes*. Springer, Berlin. MR1717346
[16] LIGGETT, T. M., SCHONMANN, R. H. and STACEY. A. M. (1997). Domination by product measures. *Ann. Appl. Probab.* **25** 71–95. MR1428500
[17] NEUHAUSER, C. (2001). Mathematical challenges in spatial ecology. *Notices Amer. Math. Soc.* **48** 1304–1314. MR1870633
[18] NEUHAUSER, C. and PACALA, S. W. (1999). An explicitly spatial version of the Lotka–Volterra model with interspecific competition. *Ann. Appl. Probab.* **9** 1226–1259. MR1728561
[19] SAWYER, S. (1976). Results for the stepping stone model for migration in population genetics. *Ann. Probab.* **4** 699–728. MR0682605
[20] WAKOLBINGER, A. (1995). Limits of spatial branching populations. *Bernoulli* **1** 171–189. MR1354460

WEIERSTRASS INSTITUTE
FOR APPLIED ANALYSIS AND STOCHASTICS
MOHRENSTR. 39
10117 BERLIN
GERMANY
E-MAIL: birkner@wias-berlin.de

TECHNISCHE UNIVERSITÄT BERLIN
INSTITUT FÜR MATHEMATIK
STRASSE DES 17. JUNI 136
10623 BERLIN
GERMANY
E-MAIL: depperschmidt@math.tu-berlin.de